\documentclass{article}
\usepackage[english]{babel}
\usepackage{graphicx} 
\usepackage{blindtext}
\usepackage{amsthm}
\usepackage{amsmath}
\usepackage{amssymb}
\usepackage{enumerate}
\usepackage{mathtools}
\usepackage{bm}
\usepackage{graphicx}
\usepackage{xcolor}
\usepackage{siunitx}
\usepackage{multirow}
\usepackage{lipsum}
\usepackage{booktabs}
\usepackage{authblk}

\title{Sections and Chapters}

\usepackage{bookmark} 
\usepackage{hyperref} 
\hypersetup{%
        pdfpagemode={UseOutlines},
        bookmarksopen,
        pdfstartview={FitH},
        colorlinks,
        linkcolor={red},
        citecolor={blue},
        urlcolor={cyan}
        }
\usepackage{float} 
\usepackage[autostyle,italian=guillemets]{csquotes}
\usepackage[backend=biber]{biblatex}
\bibliography{bibliografia.bib}

\title{Inverse modeling of porous flow through deep neural networks:\\
the case of coffee percolation}

\author[1]{A. Barletta}
\author[3]{S. Cuomo}
\author[2]{N. Egidi}
\author[2]{J. Giacomini}
\author[2]{P. Maponi}

\affil[1]{Officina Digitale Department, Spici srl, Naples, Italy\\
\texttt{barletta@spici.eu}}
\affil[2]{School of Science and Technologies - Mathematics Section, Università degli Studi di Camerino}
\affil[3]{Università degli studi di Napoli Federico II, Naples, Italy\\
\texttt{salvatore.cuomo@unina.it}}

\date{\today}

\begin{document}
\maketitle
\begin{abstract}
This work addresses the inverse problem of espresso coffee extraction, in which one aims to reconstruct the brewing conditions that generate a desired chemical profile in the final beverage. Starting from a high-fidelity multiphysics percolation model, describing fluid flow, solute transport, solid, liquid reactions, and heat exchange within the coffee bed, we derive a reduced forward operator mapping controllable brewing parameters to the concentrations of the main chemical species in the cup.

From a mathematical standpoint, we formalize the structural requirements for the local solvability of inverse problems, providing a \emph{minimal analytical condition} for the existence of a (local) inverse map: continuous differentiability of the forward operator and a locally constant, nondegenerate Jacobian rank. Under these assumptions, the Constant Rank Theorem ensures that the image of the forward operator is a smooth embedded manifold on which well-defined local right-inverses exist.

Because the forward operator is accessible only through an expensive numerical simulator, we construct a neural surrogate model that reproduces it with high accuracy ($R^2>0.99$ across all solutes). We then learn an approximate right-inverse through a composite objective enforcing parameter reconstruction and forward consistency. A tailored data-augmentation pipeline, based on mixture interpolation and spline-based temperature interpolation, densifies the sampling of the recipe space and significantly improves generalisation.

Extensive experiments, including off-grid validation, show that the learned inverse map accurately reconstructs brewing temperature, grind size, and powder composition. The resulting framework combines rigorous analytical guarantees with modern data-driven methods, providing a principled and computationally efficient solution to the inverse extraction problem and enabling personalised brewing, recipe optimisation, and integration into smart coffee-machine systems.
\end{abstract}

\section{Introduction}

Coffee is among the most widely consumed beverages globally, with a supply chain that spans nearly every country, from cultivation in climatically favorable regions to roasting, distribution, and consumption worldwide. Although coffee production is geographically constrained, its cultural and economic relevance is universal: virtually every nation includes coffee drinkers among its population, each adopting distinct extraction methods shaped by local traditions and technological preferences. The motivations behind coffee consumption are multifaceted, encompassing functional, social, and health-related dimensions. Coffee serves not only as a stimulant but also as a social catalyst and a source of physiological benefits. These diverse roles underscore the growing demand for customized coffee preparation, tailored to individual preferences that reflect cultural heritage, lifestyle, and even transient physiological states. For example, a strong brew may be favored to enhance alertness after sleep deprivation, while a lighter infusion might be preferred following a heavy meal. In this context, the coffee industry is increasingly oriented toward personalized extraction technologies, aiming to optimize sensory profiles and bioactive compound delivery in response to consumer-specific needs.

Scientific research on coffee, promoted by the interests of operators in the coffee sector, is stimulated by the multidisciplinary complexity behind the description of the extraction process. From a biological and chemical standpoint, coffee is a chemically heterogeneous system exhibiting low stability with respect to volatile compounds and foam persistence, and distinguished by a complex and distinctive sensory profile. Numerous factors, such as type of coffee blend, roasting and grinding degree and preparation method (i.e., coffee extraction technique) can influence the levels of chemical compounds that arrive in the cup and, accordingly, the final coffee flavour, quality and health-related aspects. From a mathematical standpoint, coffee extraction is a non-trivial process that can be described as the percolation of a liquid in a porous medium, formalized through Fluid Dynamics equations. Some percolation models and their analysis are shown in~\cite{giacomini2020,giacomini2023,giacomini2022cmmse,giacomini2022TDS,moroney2015modelling,moroney2016,moroney2016_3d,moroney2019,cameron2020}. 

While the \emph{forward problem}—predicting the coffee composition resulting from a given set of brewing parameters—is already challenging due to the nonlinear, multiphase nature of the process, an even more intricate task is the associated \emph{inverse problem}:
\begin{quote}
\emph{Given a desired sensory or chemical target for the extracted coffee, determine the brewing parameters that produce it.}
\end{quote}
Inverse extraction lies at the core of beverage personalization, enabling applications such as recipe optimization, adaptive brewing in smart machines, and quality control within industrial settings. However, inverse problems are typically ill-posed and highly sensitive to perturbations~\cite{Tikhonov1977, DeVito2005, Kaipio2004}, especially when the underlying forward map is non-injective, high-dimensional, or accessible only through expensive numerical simulations.

In the context of coffee percolation, the forward operator mapping brewing parameters to cup chemistry is not available in closed form but only through a high-fidelity numerical simulator derived from the multiphysics percolation model. Therefore, constructing its inverse requires careful consideration of the mathematical properties of the forward map, such as local invertibility and manifold structure, and the use of data-driven techniques to approximate stable right-inverse selections.

\paragraph{Approach to the inverse problem.}
The present work addresses the inverse problem through a hybrid analytical–data-driven methodology:

We first interpret the forward map as a smooth operator between finite-dimensional parameter and output spaces, and we verify that it satisfies the differentiability and local-constant-rank conditions required by the Constant Rank Theorem. This ensures the existence of locally well-defined right inverses on suitable output manifolds.

Because the full percolation simulator is computationally expensive and non-differentiable with respect to its inputs, we train a neural surrogate model to approximate the forward operator on a dense, augmented dataset generated through simulation and physically informed interpolation.

We then construct an approximate inverse map by training a neural network that enforces consistency with the surrogate forward model, ensuring $f(g(y))\approx y$ while respecting physical constraints and interpretability of the reconstructed brewing parameters.

The resulting framework combines theoretical guarantees of local invertibility with the flexibility of modern machine-learning techniques, enabling robust inversion even in rank-deficient or non-linear regions of the parameter space.

\paragraph{Structure of the paper.}
The remainder of the paper is organized as follows:
\begin{itemize}
\item In Section \ref{sec:background}, we recall the mathematical foundations of forward and inverse problems, describe local invertibility conditions, and discuss their empirical verification in data-driven contexts.
\item Section \ref{sec:From analytical to data-driven inversion} introduces the multiphysics percolation model, the associated numerical simulator, and the construction of the reduced forward operator. We then describe the surrogate model and provide a detailed analysis of the Constant Rank assumptions.
\item Section~\ref{sec:Experimental_results} describes the construction of the datasets, including the structured simulation grid, the off-grid complementary set, and the data augmentation pipeline, together with the training procedure for both the surrogate forward operator and the inverse model. The section also reports the full set of experimental results, assessing the performance of the surrogate and the accuracy and generalization properties of the learned inverse map.
\item Conclusions and possible future developments are summarized in Section \ref{sec:conclusions}
\end{itemize}

\section{Mathematical background}\label{sec:background}

\subsection{Forward and inverse problems}
Let $f : \mathcal{X} \to \mathcal{Y}$ denote a deterministic mapping between two metric spaces, typically representing a physical or computational model. 
The associated \emph{forward problem} consists in evaluating the response $y \in \mathcal{Y}$ corresponding to a given input or parameter vector $x \in \mathcal{X}$, according to
\[
y = f(x).
\]
Forward problems are, in most physically meaningful settings, well-posed in the sense of Hadamard (see, e.g., \cite{Tikhonov1977}), namely:  
\begin{enumerate}[a)]
    \item a solution exists,
    \item it is unique,
    \item it depends continuously on the data.
\end{enumerate}
Conversely, the \emph{inverse problem} aims at recovering the cause $x$ that produced a given observation $y$, i.e.
\begin{equation}\label{cause_effect}
    x = f^{-1}(y),
\end{equation}
or, more generally, at finding an $x$ such that $f(x) \approx y$ when the model or the measurements are noisy. 
Inverse problems are frequently ill-posed, in the sense that one or more of the Hadamard conditions fail: that is, existence, uniqueness, or stability may be lost \cite{Tikhonov1977, DeVito2005}. 
Regularization, prior information, or structural constraints are then introduced to stabilize the inversion process, either deterministically through Tikhonov-type penalization \cite{Tikhonov1977, DeVito2005}, or probabilistically within a Bayesian framework \cite{Kaipio2004, MohamadDjafari2013}.

\paragraph{Examples and context.}
Inverse problems appear in numerous areas of science and engineering, including system identification, optics, radar, acoustics, communication theory, signal processing, medical imaging, computer vision~\cite{MohamadDjafari2013,Pizlo2001}, geophysics, meteorology, astronomy, remote sensing, and machine learning~\cite{DeVito2005}. 
They are also central in applied engineering contexts such as nondestructive testing and slope stability analysis \cite{Cardenas2019}.

\subsection{Mathematical formulation}
We now introduce a general framework for the analysis of forward and inverse mappings.

Let $\Omega_{\mathrm{in}}\subseteq \mathbb{R}^n$ denote the \emph{input space} (e.g., model parameters or control variables),we consider a sufficiently smooth map
\begin{equation}\label{eq:forward_general} 
    f:\; \Omega_{\mathrm{in}} \longrightarrow \mathbb{R}^m,
\end{equation}
which defines the forward operator.  

The image of $f$,
\begin{equation}\label{eq:image_general}
  \Omega_{\mathrm{out}} \coloneqq \mathrm{Im}(f) = \big\{\, y\in\mathbb{R}^m\ \big|\ \exists\, x\in\Omega_{\mathrm{in}}\ \text{with } y=f(x)\,\big\},
\end{equation}
represents the set of \emph{attainable outputs} or physically realizable configurations.

The \emph{inverse problem} consists in recovering, for a given observation $y\in\Omega_{\mathrm{out}}$, a corresponding $x\in\Omega_{\mathrm{in}}$ satisfying $f(x)=y$.  
If $f$ is bijective, the inverse map $f^{-1}:\Omega_{\mathrm{out}}\to\Omega_{\mathrm{in}}$ exists and is single-valued.  
However, in most applications $f$ is neither injective nor surjective, and thus one introduces the \emph{set-valued inverse}
\begin{equation}\label{eq:set_valued_inverse_general}
    f^{-1}(y)\; \coloneqq\; \{\, x\in \Omega_{\mathrm{in}} \;:\; f(x)=y \,\},
    \qquad y\in \Omega_{\mathrm{out}},
\end{equation}
which may contain multiple preimages or be empty outside the attainable set.  
In such cases, one often seeks a \emph{selection function}
\begin{equation}\label{eq:right_inverse_general}
    g:\Omega_{\mathrm{out}}\to \Omega_{\mathrm{in}}
    \quad\text{such that}\quad
    f\circ g = \mathrm{Id}_{\Omega_{\mathrm{out}}},
\end{equation}
that is, a right inverse of $f$ restricted to its image.  
The choice of $g$ may encode additional regularity, optimality, or prior information, depending on the context.

\subsection{Local invertibility and the Constant Rank Theorem}
When $\Omega_{\mathrm{in}}\subseteq \mathbb{R}^n$ and $\Omega_{\mathrm{out}}\subseteq \mathbb{R}^m$, the local invertibility of $f$ depends on the relative dimensions $n$ and $m$ and on the rank of its Jacobian $Df(x)$.  

If $n=m$ and $\det Df(x^\ast)\neq 0$, the classical Inverse Function Theorem ensures that $f$ is a local diffeomorphism near $x^\ast$ (see, \cite{Spivak1965}).  

If instead $n\neq m$, the Inverse Function Theorem does not apply directly. In this more general setting, the \emph{Constant Rank Theorem} provides a suitable geometric description (see, e.g., \cite{Lee2013}).
  
Suppose $f\in C^k(\Omega_{\mathrm{in}},\Omega_{\mathrm{out}})$ for some $k\ge 1$, and that $Df(x^\ast)$ has constant rank $r\le \min\{n,m\}$ in a neighborhood of $x^\ast$.  
Then the image of $f$ near $y^\ast=f(x^\ast)$ is a $r$-dimensional embedded 
submanifold $\mathcal{M}\subset\mathbb{R}^m$, and there exist neighborhoods 
$U\subset\Omega_{\mathrm{in}}$ and $V\subset\Omega_{\mathrm{out}}$ such that
\begin{equation}\label{eq:rank_thm_general}
    f:\ U \longrightarrow V\cap\mathcal{M}
\end{equation}
is a $C^1$ diffeomorphism onto its image, i.e.\ a local diffeomorphism 
\emph{between $U$ and the submanifold $V\cap\mathcal{M}$}.
  
Consequently, a local inverse $f^{-1}:(V\cap\mathcal{M})\to U$ exists and is $C^1$.  
In geometric terms, $f$ acts locally as a smooth embedding of an $r$-dimensional manifold into $\mathbb{R}^m$ (Ref. \cite{Lee2013}).
This framework encompasses both the classical case $n=m$ (local invertibility on open sets) and the rank-deficient case $n<m$ (invertibility along lower-dimensional manifolds). Specific applications correspond to particular choices of $n$, $m$, and of the physical interpretation of $f$.  
\subsection{Modelling stance and verifiability}

Trivializing the square case $n=m$, where the classical Inverse Function Theorem guarantees local invertibility whenever $\det Df(x^\ast)\neq 0$, we focus here on the more general and practically relevant situation in which the input and output spaces have different dimensions. 
This is the typical setting in inverse problems where the number of measurable quantities exceeds or falls short of the number of controllable parameters.

In such cases, the existence of a local inverse for a mapping $f:\Omega_{\mathrm{in}}\to\Omega_{\mathrm{out}}$ requires, in its minimal form:
\begin{enumerate}[(i)]
    \item \textbf{Differentiability:} $f$ is continuously differentiable in a neighborhood of interest; 
    \item \textbf{Local constant rank:} the Jacobian $Df(x)$ maintains a locally constant, nondegenerate rank $r\le \min\{n,m\}$.
\end{enumerate}
Under these conditions, $f$ behaves locally as a diffeomorphism between an open subset of $\Omega_{\mathrm{in}}$ and an $r$-dimensional manifold $\mathcal{M}\subseteq\Omega_{\mathrm{out}}$, allowing for a local right inverse on $\mathcal{M}$.

\paragraph{From theory to practice.}
In real-world scenarios, one rarely has an explicit analytical expression of $f$. 
Depending on the context, $f$ may be accessible only:
\begin{enumerate}[a)]
    \item through a \emph{numerical simulator}, mapping inputs $x$ to computed outputs $f(x)$ with possible numerical or discretization errors;
    \item or through a \emph{collection of empirical observations} $\{(x_i,y_i)\}$, where $y_i$ are sensor measurements or experimental data.
\end{enumerate}
In both settings, $f$ effectively acts as a \emph{black-box operator}: it can be evaluated, but not differentiated analytically.

\paragraph{Empirical assessment of local invertibility.}
The mathematical hypotheses above can be examined empirically through numerical diagnostics that depend on the available information:

\begin{description}
    \item[Derivative-aware regime:] 
    If derivatives, or their approximations, are available via adjoint models, automatic differentiation, or finite differences, one may estimate the rank and conditioning of $Df(x)$ using singular value decomposition (SVD). 
    The persistence of a nonzero smallest singular value $\sigma_r(Df(x))$ in a neighborhood of $x$ supports the assumption of local constant rank and numerical stability.

    \item[Derivative-free regime:] 
    When only point evaluations of $f$ are accessible, local structure can be probed by performing a \emph{local Principal Component Analysis (LPCA)} on nearby outputs $\{y_j=f(x_j)\}$. 
    The number of dominant singular values of the local covariance matrix yields an estimate of the intrinsic dimension $r$ of the image manifold $\mathcal{M}$. 
    Alternative manifold-learning methods, such as Locally Linear Embedding (LLE) \cite{RoweisSaul2000}, Laplacian Eigenmaps \cite{BelkinNiyogi2003}, or the Local Dimensionality Estimator (LDE) \cite{LevinaBickel2004}, provide complementary tools to infer local dimensionality and validate the constant-rank hypothesis empirically.

\end{description}

\paragraph{Interpretation.}
These empirical analyses offer a practical way to assess whether the operational assumptions of differentiability and local constant rank are plausibly satisfied in a given region of the data space.  They do not replace the formal mathematical conditions, but translate them into verifiable hypotheses that can be tested before attempting the construction of an analytical, numerical, or data-driven inverse.

\section{From analytical to data-driven inversion}\label{sec:From analytical to data-driven inversion}
Based on the description of percolation that is usually applied in the field of hydrogeology, the percolation process that occurs during espresso coffee extraction can be described by a similar model. From a Fluid Dynamics standpoint, the model encompasses fluid motion, particle detachment from the porous matrix, dissolved chemical species, and heat exchange between the fluid and solid medium. The proposed model relies on the following assumptions: a) the porous medium is considered isotropic and homogeneous, b) the medium is assumed to be saturated, with no presence of a gaseous phase, and local thermal equilibrium is established between the ground coffee and the infiltrating water. Thus, the model reads:
\begin{equation}\label{perc_model}
\begin{cases}
  \displaystyle  S_0\frac{\partial h}{\partial t}+\nabla\cdot\bm{q} = 0,\\
    \bm{q}=-\bm{K}f_\mu\cdot\left( \nabla h+\chi\bm{e}\right),\\ \vspace{0.1cm}
    \displaystyle \epsilon\frac{\partial C_k}{\partial t}+\bm{q}\cdot\nabla C_k+\nabla\cdot\bm{j_k}=R_k, \qquad k=1,\dots,N_{s},\\ \vspace{0.2cm}
   \displaystyle  \epsilon_s\frac{\partial C_k^s}{\partial t}=R_k^s, \qquad k=1,\dots,N_{s},\\ 
   \displaystyle  \left(\epsilon\rho c+\epsilon_s\rho^s c^s \right)\frac{\partial T}{\partial t}+\rho c \bm{q}\cdot \nabla T-\nabla \cdot \left( \bm{\Lambda} \cdot \nabla T \right)=0,
\end{cases}
\end{equation}

where all the equations are prescribed in an open spatial domain $\mathcal{C}$ with radial symmetry, representing the coffee pod, and for $t\in(0,\tau)$, with $\tau>0$ the percolation time. In saturated porous media, fluid flow is governed by Richards’ equation~\cite{spadoni2018}, derived from the combination of the mass conservation law (first equation in~\eqref{perc_model}) and Darcy’s law (second equation in~\eqref{perc_model}), and expressed in terms of the hydraulic head $h$. The hydraulic head is defined as $h=\psi+x_3$, where $\psi$ is the pressure head and $x_3$ is the height into $\Omega$. The specific storage coefficient $S_0$ represents compressibility. The Darcy flux $q$ is described by Darcy’s law, incorporating the hydraulic conductivity tensor $\bm{K}$, the unit vector $e=(0,0,1)^T$, where $T$ stands for the transposition operator, and the coefficients $f_\mu$ and $\chi$, which model the influence of temperature and pressure on fluid viscosity and density, respectively. The third governing equation is an advective–diffusive–reactive equation describing the transport, diffusion and dissolution or erosion of the $k$-th chemical species. Here, the unknown is the mass concentration of the considered species, which may pertain to either solid or liquid phases. In these equations, $\epsilon$ is the porosity of the medium, $\bm{j_k}$ is the hydrodynamic diffusion-dispersion, defined by Fick’s law $\bm{j_k}=\bm{D_k}\cdot\nabla C_k$, where $\bm{D_k}$ is the hydrodynamic dispersion tensor \cite{giacomini2023}, consisting of the sum of molecular diffusion and mechanical diffusion, $R_k$ denotes the total reaction rate of species $k$. The fourth governing equation describes the mass balance for solid species $k$, characterized by the unknown concentration $C_k^s$. The species is described as solid because it is bound to the solid matrix, hence it is excluded from transport and diffusion processes. Here, $\epsilon_s=1-\epsilon$ is the solid volume fraction and $R_k^s$ denotes the total reaction rate of the $m$-th solid species. The total reaction rate terms $R_k,R_k^s$ are defined as:
\[
R_k=\alpha_k(1-\epsilon)C_k^s, \quad R_k^s=-\alpha_k(1-\epsilon)C_k^s,
\]
where the coefficients $\alpha_k$ are the following functions of incoming water pressure and temperature $p_{z0},T_{z0}$, respectively:
\[
\alpha_k =A_0 + aT_{z0} + bp_{z0} + cT_{z0}^2 + dp_{z0}^2 +fT_{z0}p_{z0} + lT_{z0}^2p_{z0} + mT_{z0}p_{z0}^2
\]
where the coefficients $A_0,a,b,c,d,f,l,m$ depend on the chemical species, the granulometry of coffee powder, and the coffee variety. The last equation in system~\eqref{perc_model} is the convective and diffusive heat equation. $T_0$ is a reference temperature, $\rho c$ is the fluid volumetric heat capacity, $\rho_s c_s$ is the solid volumetric heat capacity and $\Lambda$ is the thermal hydrodynamic conductivity tensor~\cite{giacomini2023}.

Let $\Omega=\mathcal{C}\cup\Gamma_1\cup\Gamma_2\cup\Gamma_3$ the closed domain, consisting of the open domain $\mathcal{C}$, the top surface $\Gamma_1$, the lateral surface $\Gamma_2$ and the bottom surface $\Gamma_3$. System~\eqref{perc_model} is endowed with the following boundary and initial conditions:
 \begin{equation}
\begin{cases}
    h=h_{z0}, & \text{on } \Gamma_1,\ t>0,\\
    \frac{\partial h}{\partial r}=0, & \text{on } \Gamma_2,\ t>0,\\
    \bm{q}\cdot \bm{n}=-\Phi_h\text{min}\{h_C-h,0\}, & \text{on } \Gamma_3,\ t>0,\\
    p=p_0(x_3), & \text{on } \mathcal{D},\ t=0,
\end{cases}
\end{equation}
\begin{equation}
\begin{cases}
    \nabla C_k\cdot \bm{n}=0, \qquad & \text{on } \Gamma_1,\Gamma_2,\ t>0,\\
    -\left(\bm{D_k}\cdot\nabla C_k\right)\cdot \bm{n}=-\Phi_k\text{min}\{C_{kC}-C_k,0\}, & \text{on } \Gamma_3, t>0,\\
    C_k=0, & \text{on } \mathcal{D},\ t=0,
\end{cases}
\end{equation}
\begin{equation}
    C_k^s=C_{k0}^s, \qquad \text{on } \mathcal{D},\ t>0,
\end{equation}
\begin{equation}\label{BIC_temp}
\begin{cases}
    T=T_{z0}, & \text{on } \Gamma_1,\ t>0,\\
    \nabla T\cdot \bm{n}=0, & \text{on } \Gamma_2,\Gamma_3,\ t>0,\\
    T=T_0, & \text{on } \mathcal{D},\ t=0,
    \end{cases}
\end{equation}
where ${\bm n}$ is the normal vector exiting a closed surface belonging to $\partial \mathcal{D}=\Gamma_1\cup\Gamma_2\cup\Gamma_3$, $r$ is the radial coordinate, $h_C,h_{z0},p_0,\Phi_h,C_{kC},\Phi_k,C_{k0}^s,T_{z0},T_{0}$ are prescribed functions or constants. Details about the derivation of system \eqref{perc_model}-\eqref{BIC_temp} can be found in \cite{giacomini2020,giacomini2023}. The system consists of dozens of parameters that depend on the geometric characteristics of the domain and the physical and chemical characteristics of the extraction process. This feature makes it a true digital twin capable of responding to the variable conditions under which the described process takes place. System \eqref{perc_model}-\eqref{BIC_temp} undergoes an approximation procedure, both in space and time. The one considered in~\cite{giacomini2020,giacomini2023} consists in a finite element approximation in space and the Adams–Bashforth/Crank–Nicolson predictor–corrector strategy in time. The final linear algebraic system is solved using a standard iterative method. 

In order to fix ideas, such a percolation model produces an \emph{in-silico} coffee as follows. 
Once the geometrical configuration of the pod is fixed, the user specifies the externally controllable brewing variables: the temperature and pressure profiles of the incoming water,  the granulometry of the powder, the mass fractions of the coffee blend, and the chemical characterization of each coffee type (i.e., the initial concentrations of all considered chemical species). Given these inputs, the numerical solver of system \eqref{perc_model}-\eqref{BIC_temp} computes the full spatio-temporal fields of pressure head, temperature, and solute concentrations inside the coffee bed, eventually yielding the chemical composition of the extracted beverage. To formalize this process, it is convenient to regard the numerical solver as an operator
\begin{equation}\label{eq:Numerical_solver_operator}
    \mathcal{S}:\;(\mathbf{p},t)\longmapsto c(t,\mathbf{p}),
\end{equation}
where $\mathbf{p}\in\mathcal{P}$ denotes the complete set of physical and chemical 
parameters appearing in the multiphysics model, and $c(t,\mathbf{p})$ is the vector of 
domain-integrated solute concentrations at time $t$.  
Only a subset of the parameters in $\mathbf{p}$ is directly controllable during brewing.  
These controllable quantities are collected into the reduced recipe vector $x\in\mathbb{R}^{n}$, for a certain $n \in \mathbb{N}$.
We can therefore consider a set $\Omega_{\mathrm{in}} \subseteq \mathbb{R}^n$ as the input space as defined in Section \ref{sec:background}, and also a parametrization map
\begin{equation}\label{eq:parametrization_map}
    E:\;\Omega_{\mathrm{in}}\longrightarrow\mathcal{P}
\end{equation}
which embeds the recipe variables into the full parameter space of the PDE model, assigning all non-controllable coefficients their nominal physical values. At the end of the extraction time interval, the cup chemistry is obtained by applying a projection operator
\begin{equation}\label{eq:projection_operator}
    \Pi:\;c(\cdot,\mathbf{p}) \longmapsto y\in\mathbb{R}^{m},
\end{equation}
which extracts the considered solute concentrations measurable in the final beverage. The resulting reduced forward operator is therefore
\begin{equation}\label{eq:general_coffee_forward_operator}
    f:\;\Omega_{\mathrm{in}}\longrightarrow \Omega_{\mathrm{out}}\subset\mathbb{R}^m,\qquad
    f(x)\coloneqq\Pi\big(\mathcal{S}(E(x),t_{\mathrm{end}})\big) = y,
\end{equation}
providing a compact and operational mapping from controllable brewing parameters to the chemical composition of the coffee cup.

\subsection{Forward model derived from numerical simulations}
The coffee percolation problem represents a concrete realization of the general framework introduced in Section \ref{sec:background}, in which the forward operator $f$ define in \eqref{eq:general_coffee_forward_operator} maps controllable brewing parameters to observable cup-chemistry outcomes.
Specifically, let $x \in \mathbb{R}^7$ denote the controllable brewing parameters,
\begin{equation}\label{eq:brewing_parameters}
    x=(x_T,x_p,x_\text{gran},x^1_\text{distr},x^2_\text{distr},x^3_\text{distr},x^4_\text{distr}),
\end{equation}
namely water temperature, pressure, grind size, and the mass fractions of four coffee components, and let $y \in \mathbb{R}^8$ denote the corresponding measured cup-chemistry vector,
\begin{equation}\label{eq:chemical_vector}
    y = (y_\mathrm{caf}, y_\mathrm{chl}, y_\mathrm{tri}, y_\mathrm{fer}, y_\mathrm{tar}, y_\mathrm{cit}, y_\mathrm{ace}, y_\mathrm{lip}),
\end{equation}
representing the concentrations of caffeine, chlorogenic acids, trigonelline, ferulic acid, tartaric acid, citric acid, acetic acid, and lipids.
The forward operator is:
\begin{equation}\label{eq:direct_solution}
    f:\;\Omega_{\mathrm{rcp}}\subseteq\mathbb{R}^7 
    \longrightarrow 
    \Omega_{\mathrm{chm}}\subset\mathbb{R}^8,
\end{equation}
with $\Omega_{\mathrm{rcp}}$ and $\Omega_{\mathrm{chm}}$ corresponding, respectively, to the input and output spaces $\Omega_{\mathrm{in}}$ and $\Omega_{\mathrm{out}}$ defined in Eq.~\eqref{eq:forward_general}.  
The image of $f$ defines the subset of \emph{physically attainable chemistries} within the output space,
\begin{equation}\label{eq:extractable_chemistry}
 \Omega_{\mathrm{chm}} \coloneqq 
 \big\{\, y \in \mathbb{R}^8 \;\big|\; 
 \exists\, x\in\Omega_{\mathrm{rcp}} \ \text{with}\ y = f(x) \,\big\},
\end{equation}
which corresponds to the attainable-output set $\mathrm{Im}(f)$ introduced in Eq. \eqref{eq:image_general} for the general case.

This operator $f$ itself is not analytical but derived from a numerical simulator $\mathcal{S}$ that solves the time-dependent multiphase percolation equations governing the extraction process.  Hence, $f$ acts as a compressed, steady-state representation of the full numerical process $\mathcal{S}$, mapping the recipe space $\Omega_{\mathrm{rcp}}$ to the attainable chemistry manifold $\Omega_{\mathrm{chm}}$.

\subsection{Verification of the Constant Rank assumptions}

For the percolation simulator defined in Section~\ref{sec:From analytical to data-driven inversion}, 
the reduced forward operator 
\[
f(x) = \Pi\!\left(\mathcal{S}(E(x), t_{\mathrm{end}})\right)
\]
inherits the regularity properties of the underlying physical model and of the projection operator $\Pi$.
Specifically, the simulator $\mathcal{S}$ computes time-dependent concentration fields that depend smoothly on the input parameters $E(x)$, owing to the differentiability of the governing partial differential equations and their constitutive relations.  
The projection $\Pi$, which extracts from these fields the final cup-chemistry vector at time $t_{\mathrm{end}}$, is itself a smooth functional.  
Assuming that the numerical solver $\mathcal{S}$ provides a sufficiently smooth 
parametric dependence on $E(x)$—as suggested by the regularity of the underlying 
continuum model and by the empirical stability of the simulated data—we may treat 
the reduced operator
\[
f = \Pi \circ \mathcal{S} \circ E
\]
as effectively $C^1$ on $\Omega_{\mathrm{in}}$ for the purpose of local 
invertibility analysis.

\vspace{0.5em}

Regarding the Jacobian rank condition, note that the dimensionality of the input and output spaces is $n=7$ and $m=8$, respectively.  
Since $n < m$, the map $f$ cannot be globally invertible onto an open subset of $\mathbb{R}^8$. 
However, if the Jacobian $Df(x)$ has locally constant rank $r$ in a neighbourhood of a reference point $x^\ast$ within the region of interest, then by the Constant Rank Theorem, the image of $f$ in a neighbourhood of $y^\ast = f(x^\ast)$ forms a smooth embedded submanifold $\mathcal{M}_{y^\ast} \subset \mathbb{R}^8$. 
In such a case, a local right-inverse exists on $\mathcal{M}_{y^\ast}$, and one can construct a piecewise-defined global inverse by suitably patching together local inverses across neighbouring regions.

A qualitative assessment of this assumption was conducted on the simulated dataset described in Section~\ref{sec:Experimental_results}.
Local dimensionality estimates obtained through LPCA and manifold-based diagnostics exhibited a stable behaviour across the sampled region, consistently indicating that the image of the forward operator lies on a smooth low-dimensional manifold. These observations are in agreement with the constant-rank hypothesis and support the applicability of the Constant Rank Theorem within the operational domain considered.




\subsection{Data-driven construction of the inverse map}

Given that the forward operator $f$ satisfies the smoothness and local rank conditions discussed in Section~\ref{sec:background}, 
there exists a local inverse selection
\[
g:\Omega_{\mathrm{chm}}\to\Omega_{\mathrm{rcp}},
\]
which can be defined locally on the attainable-chemistry manifold $\Omega_{\mathrm{chm}}$, as introduced in Eq. \eqref{eq:extractable_chemistry}.
Our goal is to approximate this right-inverse through a learnable map $g_\theta$ satisfying
\begin{equation}\label{eq:direct_composed_inverse}
    f\big(g_\theta(y)\big) \approx y, 
\qquad y \in \Omega_{\mathrm{chm}},
\end{equation}
while ensuring that $g_\theta(y) \in \Omega_{\mathrm{rcp}}$, i.e., that predicted brewing configurations correspond to realistic process parameters.

In this work, the inverse map is realized as a neural network $g_\theta$ trained on paired samples $\{(x_i, y_i)\}_{i=1}^{N_0}$, where each pair is obtained by evaluating the forward operator $f$ on a given input configuration, i.e., $y_i = f(x_i)$. 
Thus, the training dataset directly reflects the input-output relationship induced by the numerical percolation model.
Training is driven by a composite loss function that combines parameter reconstruction and forward-model consistency:
\begin{equation}\label{eq:loss_function_general}
    \min_{\theta}\;
    \mathcal{L}(\theta)
    =\alpha\,\mathcal{L}_{\text{par}}(\theta)
    +\beta\,\mathcal{L}_{\text{cons}}(\theta),
\end{equation}
where $\alpha,\beta>0$ are balancing coefficients.  

The first term $\mathcal{L}_{\text{par}}$ enforces proximity between the predicted and reference recipes, acting as a regularizer in the parameter space.  
The second term $\mathcal{L}_{\text{cons}}$ enforces forward consistency by penalizing discrepancies between the measured chemistry $y$ and the one predicted when the reconstructed recipe is propagated through a the forward operator $f$:
\begin{equation}\label{eq:loss_consistency}
   \mathcal{L}_{\text{cons}}(\theta)
   =\big\|f\big(g_\theta(y)\big) - y\big\|_2^2.
\end{equation}
This term enforces a \emph{cycle-consistency constraint} ensuring that an observation $y$ can be mapped to a plausible cause $x=g_\theta(y)$ and back through the forward process.  
In regions where $f$ is locally non-injective, the consistency term also acts as a selection mechanism, promoting reconstructions that are locally stable and physically coherent under the composition $f\circ g_\theta$.

Since the full simulator $S$ is computationally expensive, we replace $f$ in Eq. \ref{eq:loss_consistency} with a differentiable surrogate model $\hat{f}_\phi$, implemented as a compact neural network trained to emulate the mapping $f:\Omega_{\mathrm{rcp}}\to\Omega_{\mathrm{chm}}$.  

The forward-consistency term becomes
\begin{equation}\label{eq:surrogate_loss}
   \mathcal{L}_{\text{cons}}(\theta)
   =\big\|\hat{f}_\phi(g_\theta(y)) - y\big\|_2^2.
\end{equation}
The surrogate $\hat{f}_\phi$ provides a differentiable approximation of the forward operator $f$, enabling gradient-based optimization of $g_\theta$ while preserving the physical structure of the percolation process.
This surrogate-based, data-driven inversion framework enables the construction of a stable approximate right-inverse of $f$ on the attainable manifold $\Omega_\text{chm}$, even in the absence of an explicit analytical expression for the forward model.

The use of a composite objective combining reconstruction and consistency terms follows a well-established paradigm in learned inverse problems and physics-informed learning \cite{Adler2018, Arridge2019}. 
On compact domains, feedforward neural networks with nonlinear activation functions are universal approximators of continuous mappings (Ref. \cite{Cybenko1989, Hornik1989}), and can thus represent inverse functions when they exist and are sufficiently well-conditioned.  
Regularization mechanisms—both architectural (e.g., spectral normalization, residual connections) and loss-based (e.g., consistency constraints, prior terms)—help mitigate the intrinsic instability of inverse mappings \cite{Kutyniok2022, Arridge2019}, 
providing a data-driven analogue of classical deterministic and Bayesian regularization schemes \cite{Tikhonov1977, DeVito2005, Kaipio2004}.


\section{Experimental Results}\label{sec:Experimental_results}

This section presents the experimental evaluation of the proposed inverse-learning framework. 
We first describe the dataset and its augmentation process, then outline the training procedure and model architectures, and finally discuss the obtained results.

\subsection{Dataset Description}\label{sec:dataset_description}

The dataset employed in this study consists of a single text file containing paired measurements of brewing conditions $x$ and the corresponding cup-chemistry profiles $y$, as defined in Eqs.~\ref{eq:brewing_parameters} and \ref{eq:chemical_vector}. Each row represents a distinct numerical simulation of the coffee percolation process generated under specified physical and compositional parameters.

For clarity, we recall the explicit structure of the input vector
\[
 x=(x_T,x_p,x_\text{gran},x^1_\text{distr},x^2_\text{distr},x^3_\text{distr},x^4_\text{distr}),
\]
and of the measured cup-chemistry vector
\[
 y = (y_\mathrm{caf}, y_\mathrm{chl}, y_\mathrm{tri}, y_\mathrm{fer}, y_\mathrm{tar}, y_\mathrm{cit}, y_\mathrm{ace}, y_\mathrm{lip}).
\]
All data were generated through numerical simulations of the percolation dynamics described in Section~\ref{sec:From analytical to data-driven inversion}, systematically varying temperature, granulometry, and powder composition. The simulation campaign was designed as a regular grid over the parameter space. Temperatures were sampled at seven (approximately) equispaced values between $88$ and $98~\si{^\circ C}$, while pressure was fixed at $9~\si{bar}$ to reflect the operating constraints of the reference espresso machine, whose pressure setting is not adjustable. Granulometry was considered at three discrete levels—coarse (\texttt{G}), optimal (\texttt{O}), and fine (\texttt{F})—each corresponding to a specific numerical value. Powder compositions were enumerated over the set of all feasible mixtures of the four coffee fractions, sampled with a step of $1/6$ for each component, subject to the normalization constraint (fractions summing to unity).
Prior to training, the raw text data were parsed into a structured tabular format. The categorical variable corresponding to granulometry (\texttt{G, O, F}) was numerically encoded as
\[
\texttt{G} \mapsto 0, \quad 
\texttt{O} \mapsto 1, \quad 
\texttt{F} \mapsto 2,
\]
yielding a compact numerical representation  of each recipe $x \in \mathbb{R}^7$, paired with the corresponding chemical extraction $y \in \mathbb{R}^8$.
This design results in a structured grid spanning the variables (temperature, granulometry, composition), where each grid point is associated with a simulated cup-chemistry vector $y$ obtained through the numerical simulator. The constant-pressure assumption facilitates direct comparison with experimental extractions performed on the physical espresso machine, which were subsequently used in a dedicated sensorial validation session aimed at qualitatively assessing the plausibility and consistency of the simulated chemistry profiles.

\subsection{Complementary Off-Grid Dataset}

A complementary dataset comprising \num{153} additional simulation points was constructed to provide an independent yet structured extension of the main grid. While it involves the same experimental variables as the primary dataset, its sampling strategy was deliberately shifted so as to explore regions of the parameter space not covered by the original grid, thereby serving as an \emph{off-grid} validation domain. Although structured, this set is offset with respect to the main grid and thus probes intermediate regions of the parameter space.

Specifically, three distinct brewing temperatures---$89$, $92.5$, and $95~\si{^\circ C}$---were selected, interleaved with but not coincident with the seven grid values of the main dataset. The extraction pressure was kept fixed at $9~\si{bar}$ for consistency, and granulometry was again considered at the three categorical levels: coarse (\texttt{G}), optimal (\texttt{O}), and fine (\texttt{F}).

Powder compositions in this complementary dataset followed a simplified yet systematic scheme based on a limited set of mixture archetypes:
\begin{itemize}
\item pure compositions consisting of $100\%$ of a single coffee fraction;
\item binary mixtures composed of $75\%$ of one fraction and $25\%$ of another;
\item a uniform mixture with $25\%$ of each of the four fractions (\texttt{A}, \texttt{R}, \texttt{L}, \texttt{E}).
\end{itemize}

The Cartesian combination of temperature, granulometry, and composition settings yields a compact yet informative collection of off-grid points, each associated with a simulated cup-chemistry vector $y = f(x)$. This complementary dataset was primarily employed to assess the generalization capability of the learning models under slightly perturbed physical and compositional conditions. Furthermore, it served as an independent benchmark to evaluate the effectiveness of the data-augmentation strategy described in the following subsection.

\subsection{Data Augmentation}

Since the original dataset was constructed as a structured grid over temperature, granulometry, and powder composition, its discrete nature led to a limited coverage of the continuous input space. 
Preliminary experiments showed that the learning model, despite being formulated as a regressors on temperature and powder composition, tended to behave as a classifiers when trained solely on grid-aligned samples, resulting in poor generalization on off-grid points. To address this limitation, we designed an extensive data-augmentation pipeline aimed at densifying the sampling and smoothing the input–output manifold.

Let 
\begin{equation}\label{eq:original_dataset}
    \mathcal{D}_0 = \{(x_i, y_i)\}_{i=1}^{N_0}
\end{equation}
denote the original dataset, where $x_i$ represents a specific brewing condition and $y_i$ the corresponding chemical profile. 
The augmentation strategy consisted of two complementary procedures forming a multidimensional interpolation lattice across the parameter space. 
\paragraph{Composition-based interpolation.}
For each fixed combination of temperature, pressure, and granulometry, synthetic coffee mixtures were generated by convexly combining the chemical profiles of the pure fractions \texttt{A}, \texttt{R}, \texttt{L}, and \texttt{E}. This approach was motivated by the empirical observation that blending the extracted beverages of pure fractions in prescribed proportions yields a cup chemistry nearly indistinguishable from that obtained by extracting the corresponding mixed powders under identical brewing conditions. Formally, this assumption is expressed as
\begin{equation}\label{eq:convex_combination_coffee}
    f\!\left(\sum_{i=1}^{4} c_i\,P_i\right) \;\approx\; \sum_{i=1}^{4} c_i\,f(P_i),
\end{equation}
where $f(\cdot)$ denotes the forward operator defined in Eq.~\ref{eq:direct_solution}, $P_i$ are the pure fractions, and $c_i$ are non-negative mixture coefficients satisfying $\displaystyle\sum_i c_i = 1$.

This equivalence was quantitatively verified using the available data: for four-component blends, the mean absolute deviation across all chemical species was below $10^{-3}$, while for binary mixtures it was of the order of machine precision. These results supported the validity of modeling the cup chemistry of arbitrary blends as convex combinations of the pure-fraction profiles, providing a physically consistent foundation for the first stage of the augmentation.

Mixing coefficients were initially sampled in large number from a Dirichlet distribution with concentration parameter $\alpha = 1$, providing uniform coverage of the simplex $\mathcal{C}_{\mathrm{Dir}}\subset [0,1]^4$. 
Although the Dirichlet$(1,\ldots,1)$ distribution is uniform over the simplex, the probability of sampling extreme compositions, i.e., points close to the vertices, is very small, leading to an empirical prevalence of moderately balanced mixtures. Consequently, randomly extracting a fixed number $M$ of samples from this pool would have produced an overrepresentation of near-uniform compositions. To ensure adequate coverage of high-concentration regimes, we introduced an additional deterministic selection step.

Given the desired number $M$ of augmented samples, we defined, for each component $i=1,\dots,4$, a set of $K = M/4$ equally spaced target concentrations in the interval $[0,1]$:
\begin{equation}\label{eq:linspace_simplex_revised}
L_i = {\ell_{i,1}, \ldots, \ell_{i,K}} \subset [0,1].
\end{equation}
These sets prescribe uniformly distributed marginal concentration levels for each coordinate direction.

For every $\ell_{i,k} \in L_i$, we selected from $\mathcal{C}_{\mathrm{Dir}}$ the vector whose $i$-th component is closest to the prescribed level:
\begin{equation}\label{eq: sampled_dirichlet_revised}
    c^{(i,k)}
    = \arg\min_{c \in \mathcal{C}_{\mathrm{Dir}}} \bigl| c_i - \ell_{i,k} \bigr| .
\end{equation}
Selection was performed by an iterative cyclic procedure over the indices $i=1,\dots,4$.

At each iteration and for each $i$, one sample $c^{(i,k)}$ was selected according to Eq. \eqref{eq: sampled_dirichlet_revised} and immediately removed from the pool $\mathcal{C}_{\mathrm{Dir}}$.
The process continued until all $M$ augmented samples were collected. This recursive cyclic strategy prevents duplicate selections and guarantees a balanced exploration of all coordinate directions.

The resulting mixture dataset is 
\begin{equation}\label{eq:sampled_points}
    \mathcal{C}_{\mathrm{aug}}= \{\, c^{(i,k)} : i = 1,\ldots,4,\; k = 1,\ldots,K \, \} \subset [0,1]^4,
\end{equation}
which provides a near-uniform tiling of the simplex when projected along each axis.
More precisely, for every canonical basis vector $e_i$, the marginal projection map
\begin{equation}\label{eq:marginal_projection_dirichlet}
    \pi_i : \mathcal{C}_{\mathrm{aug}} \to [0,1], 
\qquad \pi_i(c) = c_i,
\end{equation}
is approximately uniformly distributed over $[0,1]$.  
Hence, the augmented compositions achieve a quasi-uniform marginal distribution along each axis 
of the simplex, while remaining valid convex combinations.

Finally, the augmented mixture dataset is:

\begin{equation}\label{eq:mixture_dataset}
    \mathcal{D}_{\mathrm{mix}}
    =
    \Bigl\{
        \bigl((x_T,\, x_p,\, x_{\mathrm{gran}},\, \mathbf{c}'),\;
        y' = \textstyle\sum_{i=1}^{4} c'_i\, y_i^{(\mathrm{pure})}\bigr)
        \;\Big|\;
        \mathbf{c}' \in \mathcal{C}_{\mathrm{aug}}
    \Bigr\}.
\end{equation}

\paragraph{Temperature-based augmentation.}
For each fixed configuration of pressure, granulometry, and powder composition, 
the available chemical outputs are observed at a discrete set of temperatures
\[
T_1 < T_2 < \cdots < T_m .
\]
To densify the sampling along the thermal axis, each chemical variable 
$y^{(j)}$ is approximated by a cubic spline
\[
s_j : [T_{\min}, T_{\max}] \to \mathbb{R},
\qquad j = 1,\ldots, 8,
\]
fitted to the original measurements $\displaystyle\{(x_{T,\ell}, y^{(j)}_\ell)\}_{\ell=1}^{N_0}$.

New synthetic temperatures are drawn from a uniform distribution over the admissible range,
\[
T_k \sim \mathcal{U}(T_{\min}, T_{\max}),
\]
and the corresponding chemical profile is obtained by evaluating the spline vector,
\begin{equation}\label{eq:spline_temperature_y}
    y_{\text{spline}}(T_k) = \bigl(s_1(T_k),\, s_2(T_k),\, \ldots,\, s_8(T_k)\bigr).
\end{equation}

The temperature-augmented set is therefore given by
\begin{equation}\label{eq:augmented_temperature}
    \mathcal{D}_{\mathrm{temp}}
    =
    \Bigl\{
        \bigl((T_k,\, x_{p},\, x_{\mathrm{gran}},\, x_{\mathrm{distr}}^1,\, x_{\mathrm{distr}}^2,\, x_{\mathrm{distr}}^3,\, x_{\mathrm{distr}}^4),\; y_{\text{spline}}(T_k)\bigr)
        \;\Big|\;
        T_k \sim \mathcal{U}(T_{\min}, T_{\max})
    \Bigr\}.
\end{equation}
which provides smooth intermediate trajectories between simulated temperature slices 
and increases the resolution of the dataset along the thermal dimension.

The final augmented dataset is the union of the original dataset (Eq. \ref{eq:original_dataset}), the additional mixture datset (Eq. \ref{eq:mixture_dataset}) and the new temperature dataset (Eq.\ref{eq:augmented_temperature}): 
\begin{equation}\label{eq:final_dataset}
    \mathcal{D} = \mathcal{D}_0 \cup \mathcal{D}_\mathrm{mix} \cup \mathcal{D}_\mathrm{temp}.
\end{equation}
Both procedures were implemented in Python within a custom class. In total, 3000 additional samples were generated from mixture interpolation and another 3000 from temperature interpolation, expanding the original 1759-sample dataset to 7759 points. This enriched dataset provided a denser and more continuous representation of the brewing parameter space, resulting in a marked improvement in model generalization beyond the discrete experimental grid. The improvement in interpolation performance was quantitatively confirmed using the complementary off-grid dataset, which served as an external validation set for the augmented model.
Finally, the full dataset was randomly partitioned into training ($70\%$), validation ($15\%$), and test ($15\%$) subsets, normalized via a Min--Max scaling procedure, and subsequently used to train the neural network described in the following subsection.

\subsection{Learning process}
As discussed in Section~\ref{sec:background}, the consistency term 
$\mathcal{L}_{\text{cons}}(\theta)$ in the inverse-learning formulation 
(Eq.~\ref{eq:loss_consistency}) requires a forward operator that maps brewing 
parameters to cup chemistry.  
In our setting, however, the reduced forward map 
$f = \Pi \circ \mathcal{S} \circ E$ is computationally expensive, since 
$\mathcal{S}$ involves the full numerical integration of the percolation 
equations and must be executed for each input configuration.  
This computational burden makes the direct use of $f$ inside the training loop 
infeasible, requiring instead a fast and differentiable approximation.

To this end, two learnable operators are introduced.  
A surrogate model $\hat f_\phi$ is first trained to approximate the reduced 
forward operator $f : \Omega_{\mathrm{rcp}} \to \Omega_{\mathrm{chm}}$, providing 
a computationally efficient proxy for the percolation dynamics.  
Building on this surrogate, an inverse map 
$g_\theta : \Omega_{\mathrm{chm}} \to \Omega_{\mathrm{rcp}}$ is then learned by 
combining parameter-reconstruction objectives with forward-consistency 
constraints, where the latter are evaluated using $\hat f_\phi$.

\subsubsection{Surrogate model via neural network}

The surrogate model, denoted by $\hat{f}_\phi$, provides the differentiable forward operator required for the consistency term in Eq.~\ref{eq:loss_consistency}. 
It approximates the numerical operator
\[
f:\; x \mapsto y,
\]
where $x$ represents the usual brewing parameters and $y$ the corresponding vector of extracted chemical concentrations. 
In essence, $\hat{f}_\phi$ offers a fast and differentiable proxy of the direct percolation dynamics, enabling gradient-based optimization of the inverse model $g_\theta$ through the forward-consistency loss.

To construct $\hat{f}_\phi$, a compact feed-forward neural network was trained on the main dataset described in Section~\ref{sec:dataset_description}. 
The model, as said in Section \ref{sec:From analytical to data-driven inversion} was optimized to reproduce the final-state chemical composition corresponding to each set of input parameters predicted by the forward operator $f$. 
Training was carried out in \texttt{PyTorch} using a hybrid mean-squared and mean-absolute error objective, adaptive learning-rate scheduling, and early stopping. 
Regularization techniques, including batch normalization and dropout, were employed to promote smoothness and prevent overfitting, ensuring that the surrogate captured the underlying physical trends.

The resulting network accurately reproduced the outputs of the direct simulator 
across all eight solutes, exhibiting stable generalization over the test domain.  
A quantitative assessment of its predictive accuracy, including detailed error 
metrics for each chemical component, is reported in Section~\ref{subsec:Network_metrics}.
Once validated, $\hat{f}_\phi$ was integrated as the differentiable forward operator within the inverse-learning pipeline, thereby linking the data-driven inverse map $g_\theta$ to the physical extraction process in a computationally tractable and fully differentiable manner.

\subsubsection{Inverse model}

The inverse mapping $g_\theta:\Omega_{\text{chm}}\to\Omega_{\mathrm{rcp}},$ was implemented as a compact feed-forward neural network with a shared feature backbone and four output heads corresponding to the blend-composition fractions, final brew temperature, percolation pressure, and grind-size class.

The shared backbone consisted of three fully-connected layers (256–128–64 units) with \texttt{ReLU} activations, batch normalization, and dropout regularization ($p=0.15$). 
The distribution head employed a softmax-normalized output to predict relative solute proportions, while the temperature and pressure heads performed scalar regression with sigmoid activation. 
A categorical classification head was used for predicting the granulometric class. 
This multi-task architecture encourages the network to learn coupled representations across correlated physical variables.

Training was performed in \texttt{PyTorch} using the Adam optimizer (learning rate $5\times10^{-4}$) with a \texttt{ReduceLROnPlateau} scheduler (patience = 150, factor = 0.98). 
The batch size was automatically chosen as the nearest power of two to the square root of the dataset size, balancing convergence stability and computational efficiency. 
Early stopping (patience = 300, $\Delta_{\min}=10^{-4}$) was used to prevent overfitting.

Each output head was trained with a task-specific loss, and the task losses were combined using a set of \emph{fixed} scalar weights:
\begin{equation}\label{eq:loss_multi}
    \mathcal{L}_{\text{multi}} = \sum_{i} w_i\,\mathcal{L}_i,
\end{equation}
where $w_i$ are non-negative constants chosen \textit{a priori} through validation and kept fixed throughout training. 
Preliminary experiments also tested dynamic weighting strategies (e.g., uncertainty-based and gradient-balancing schemes), but these approaches led to unstable optimization dynamics. 
In particular, some tasks were assigned excessively low weights during training, resulting in underrepresentation in the total loss and a subsequent degradation of their predictive metrics. 
For this reason, a fixed-weight configuration was adopted, which provided more stable convergence and consistent performance across all output tasks (solute distribution, temperature, pressure, and granulometry).

Several adaptive weighting strategies were explored during the preliminary phase. 
These included uncertainty-based formulations, where task variances were learned jointly with the network parameters (Ref.\cite{Kendall2018multi}), and softmax-normalized weighting schemes that adjusted the relative task importance during training.  
Although these methods occasionally accelerated convergence on dominant tasks, they produced highly unbalanced weighting distributions, causing certain objectives (especially temperature or powder distribution) to contribute negligibly to the total loss. 
This imbalance resulted in degraded performance for those targets and reduced overall network generalization. 
Consequently, the fixed-weight formulation was retained as a more stable and interpretable solution, ensuring a uniform contribution of all physical targets to the total training objective.

\paragraph{Loss terms.}  
Throughout the paper, a hat symbol denotes network predictions

The solute-distribution head was trained using a hybrid mean-squared and mean-absolute error:
\begin{equation}\label{eq:loss_distrib}
    \mathcal{L}_{\text{distrib}} = 0.8\,\|y_{\text{distrib}}-\hat{y}_{\text{distrib}}\|_2^2
+ 0.2\,\|y_{\text{distrib}}-\hat{y}_{\text{distrib}}\|_1
+ \lambda_s\|\hat{y}_{\text{distrib}}\|_1,
\end{equation}
where $\lambda_s=10^{-3}$ enforces sparsity and smoothness of the predicted distributions.

The temperature head employed an asymmetric regression loss that penalized overestimation errors more strongly, since from a sensorial perspective such errors have a markedly stronger impact on taste:
\begin{equation}\label{eq:loss_temp}
    \mathcal{L}_{\text{temp}}
    =
    0.7\, w(x)\,\|y - \hat{y}\|_2^2
    \;+\;
    0.3\, w(x)\,\|y - \hat{y}\|_1,
    \qquad
    w(x)=
    \begin{cases}
        \alpha, & \hat{y}>y,\\
        1, & \text{otherwise},
    \end{cases}
\end{equation}

with $\alpha=1.5$.

The pressure regression head used a standard mean-squared error:
\begin{equation}\label{eq:loss_press}
    \mathcal{L}_{\text{press}} = \|y_{\text{press}}-\hat{y}_{\text{press}}\|_2^2.
\end{equation}

The grind-size classification head was optimized using the categorical cross-entropy:
\begin{equation}\label{eq:loss_gran}
   \mathcal{L}_{\text{gran}} = 
    -\sum_{k=1}^{K} y_{\text{gran},k}\log \hat{y}_{\text{gran},k}. 
\end{equation}
\paragraph{Physics-based reconstruction constraint.}  
In addition to these task losses, a physical consistency term was imposed by coupling the inverse model with the direct surrogate model with the "frozen" pretrained forward model $\hat f_\phi$ 
\begin{equation}\label{eq:loss_reconstruction}
   \mathcal{L}_{\text{recon}} = 
\|\hat f_\phi(g_\theta(y)) - y\|_2^2.
\end{equation}
This term enforces that the predicted brewing parameters, when re-propagated 
through the surrogate forward model, reproduce the original chemical profile $y$, 
thereby regularizing the solution space and promoting physically consistent inversions.

The final objective combined all contributions as:
\begin{equation} \label{eq:total_loss_inverse}
    \mathcal{L}_{\text{total}} =
    \mathcal{L}_{\text{multi}} + \beta\,\mathcal{L}_{\text{recon}}, \quad \beta=3.
\end{equation}

Training typically converged after $3$–$4\times10^3$ epochs, after which the best-performing model (minimum validation loss) was retained.

\subsection{Network metrics}
\label{subsec:Network_metrics}

This section reports the quantitative evaluation of both the surrogate forward model 
$\hat f_\phi$ and the inverse model $g_\theta$, using the test sets described in 
Section~\ref{sec:Experimental_results}. All figures included in this section are directly 
generated from the numerical experiments and are intended to illustrate predictive accuracy 
across the various chemical species and brewing parameters.

\subsubsection{Forward surrogate model performance}

The surrogate model $\hat f_\phi$ was evaluated on the held-out test dataset. 
Global metrics as well as per–solute metrics are reported in Table~\ref{tab:forward_metrics}. 
The results show high predictive accuracy across all eight solutes, with $R^2$ values 
consistently above $0.99$.

\begin{table}[H]
\centering
\caption{Performance metrics for the surrogate forward model.}
\label{tab:forward_metrics}
\begin{tabular}{lccc}
\toprule
\textbf{Species} & \textbf{MSE} & \textbf{MAE} & \textbf{$R^2$} \\
\midrule
Caffeine              & 3.0e-06   & 1.334e-03 & 0.9989 \\
Chlorogenic Acids     & 5.0e-06   & 1.773e-03 & 0.9985 \\
Trigonelline          & 3.0e-08   & 4.20e-04  & 0.9984 \\
Citric Acid           & 1.04e-04  & 7.769e-03 & 0.9979 \\
Acetic Acid           & 3.10e-05  & 4.490e-03 & 0.9908 \\
Tartaric Acid         & 2.3e-08   & 3.57e-04  & 0.9986 \\
Ferulic Acid          & 6.1e-09   & 5.8e-05   & 0.9990 \\
Lipids                & 7.0e-06   & 2.164e-03 & 0.9965 \\
\bottomrule
\end{tabular}
\end{table}

Representative ``predicted vs.true'' plots for each chemical species are shown below.  

\begin{figure}[H]
\centering
\includegraphics[width=0.49\textwidth]{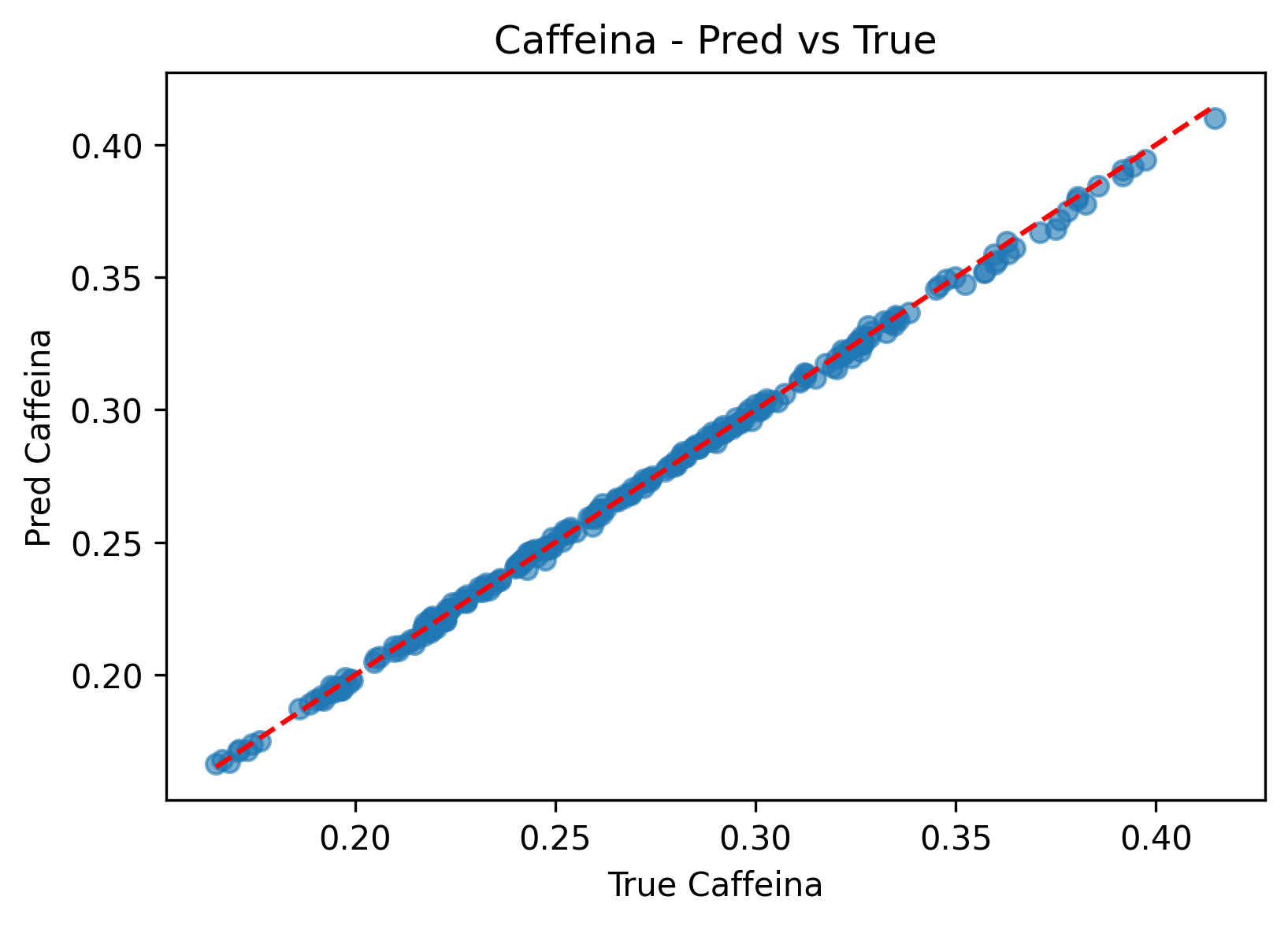}
\includegraphics[width=0.49\textwidth]{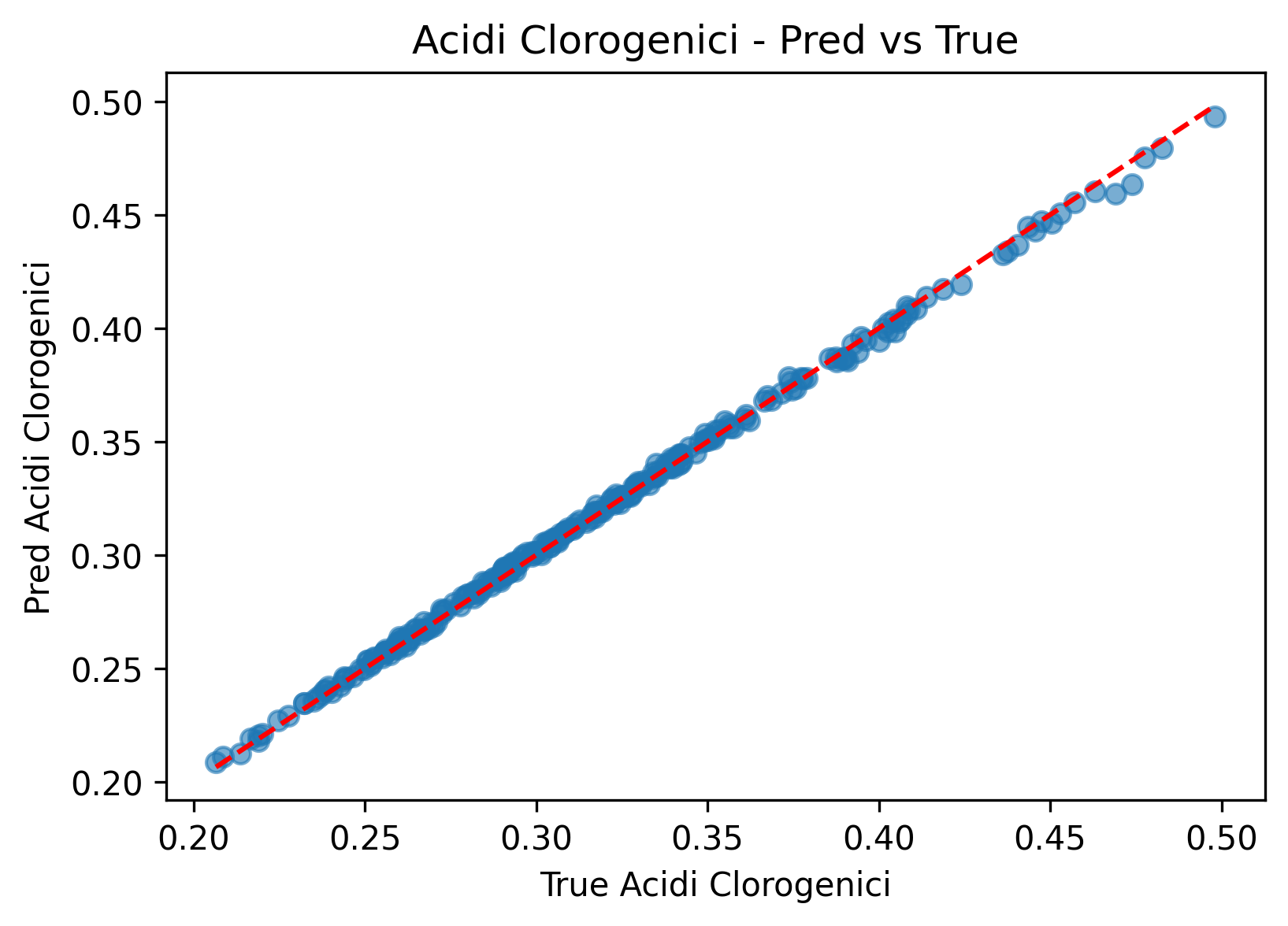}
\caption{Surrogate forward model: predicted vs.\ true concentrations for Caffeine and Chlorogenic Acids.}
\end{figure}

\begin{figure}[H]
\centering
\includegraphics[width=0.49\textwidth]{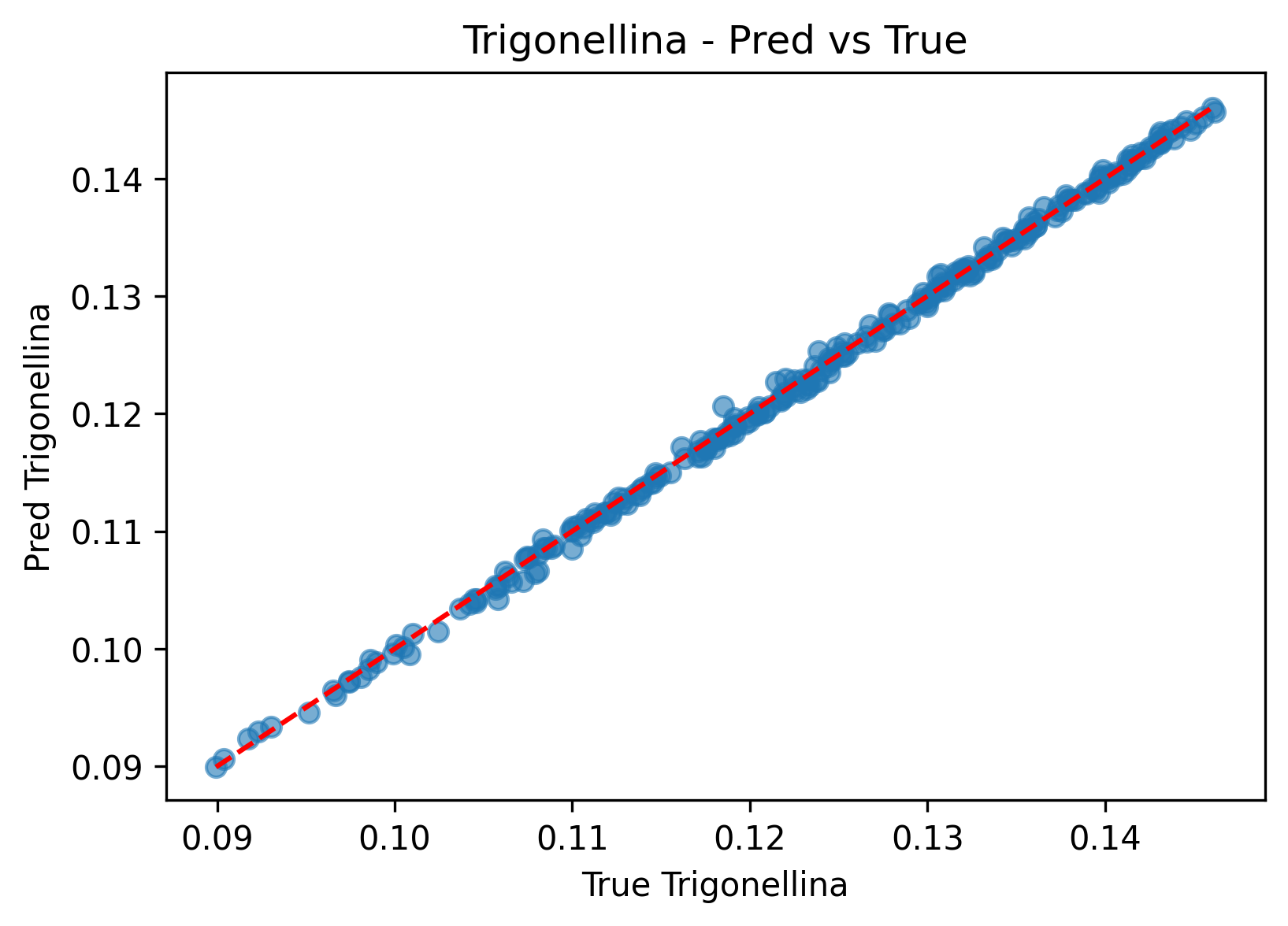}
\includegraphics[width=0.49\textwidth]{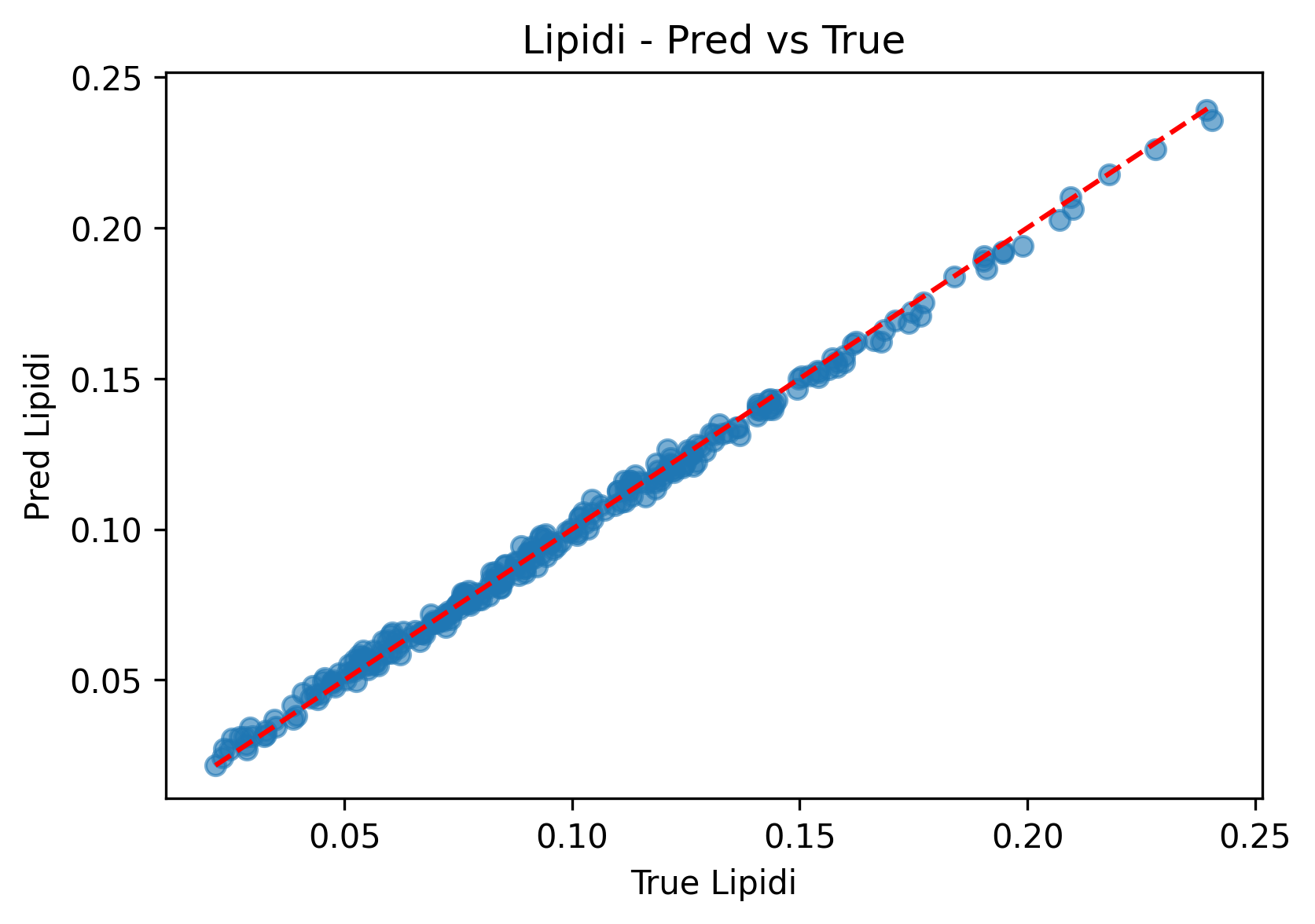}
\caption{Surrogate forward model: predicted vs. true concentrations for Trigonelline and Lipids.}
\end{figure}

\begin{figure}[H]
\centering
\includegraphics[width=0.49\textwidth]{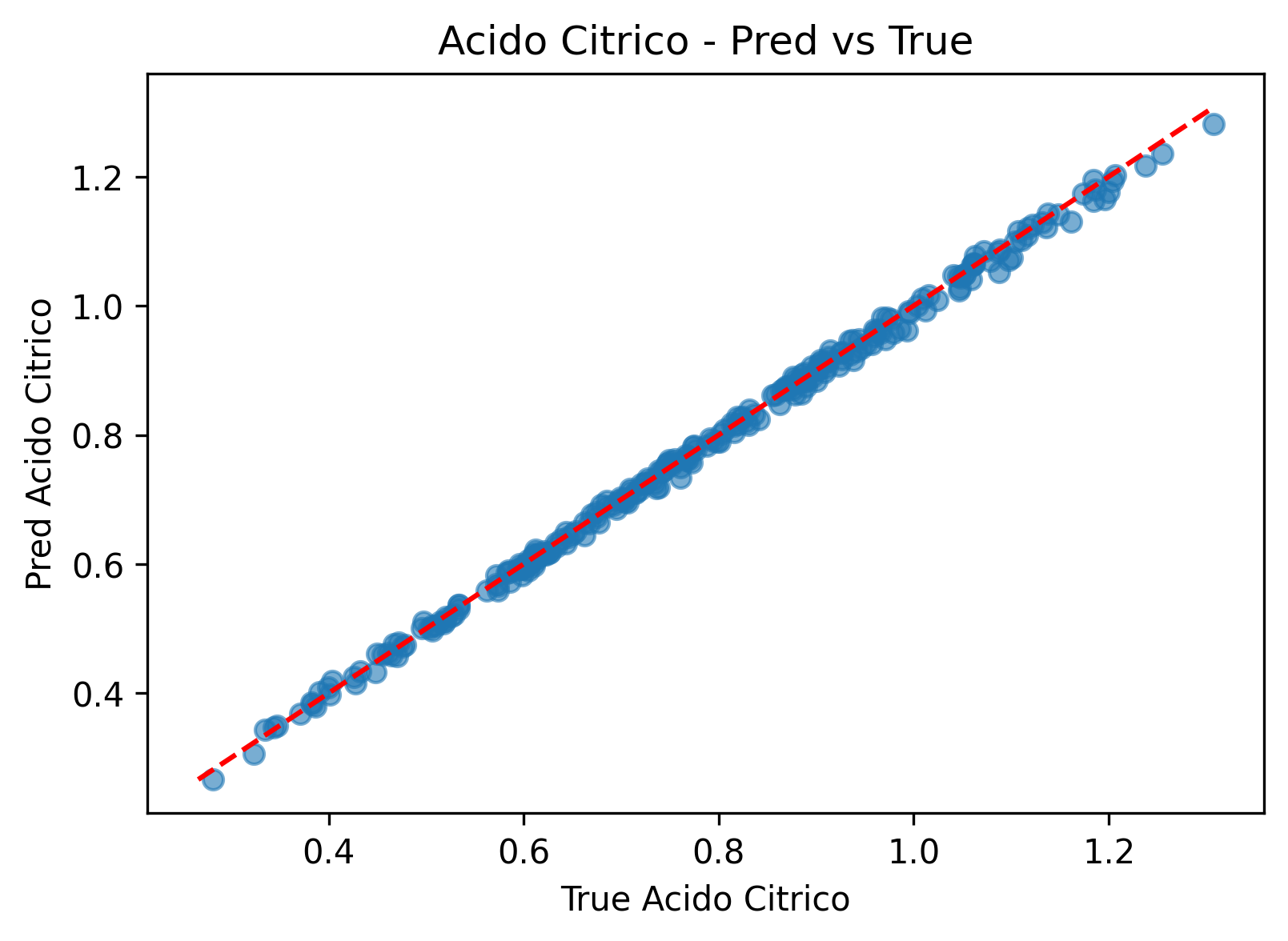}
\includegraphics[width=0.49\textwidth]{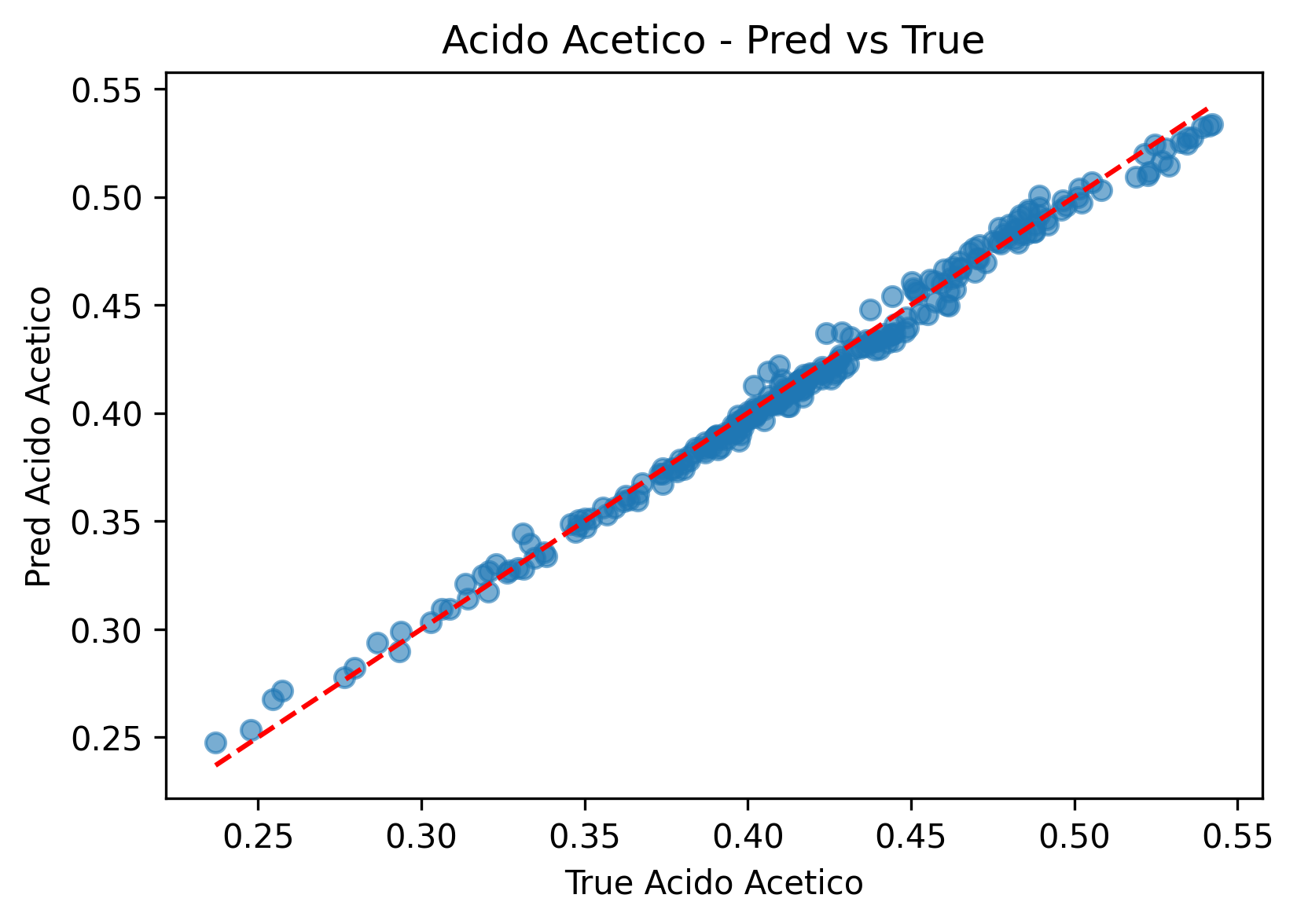}
\caption{Surrogate forward model: predicted vs. true concentrations for Citric Acid and Acetic Acid.}
\end{figure}

\begin{figure}[H]
\centering
\includegraphics[width=0.49\textwidth]{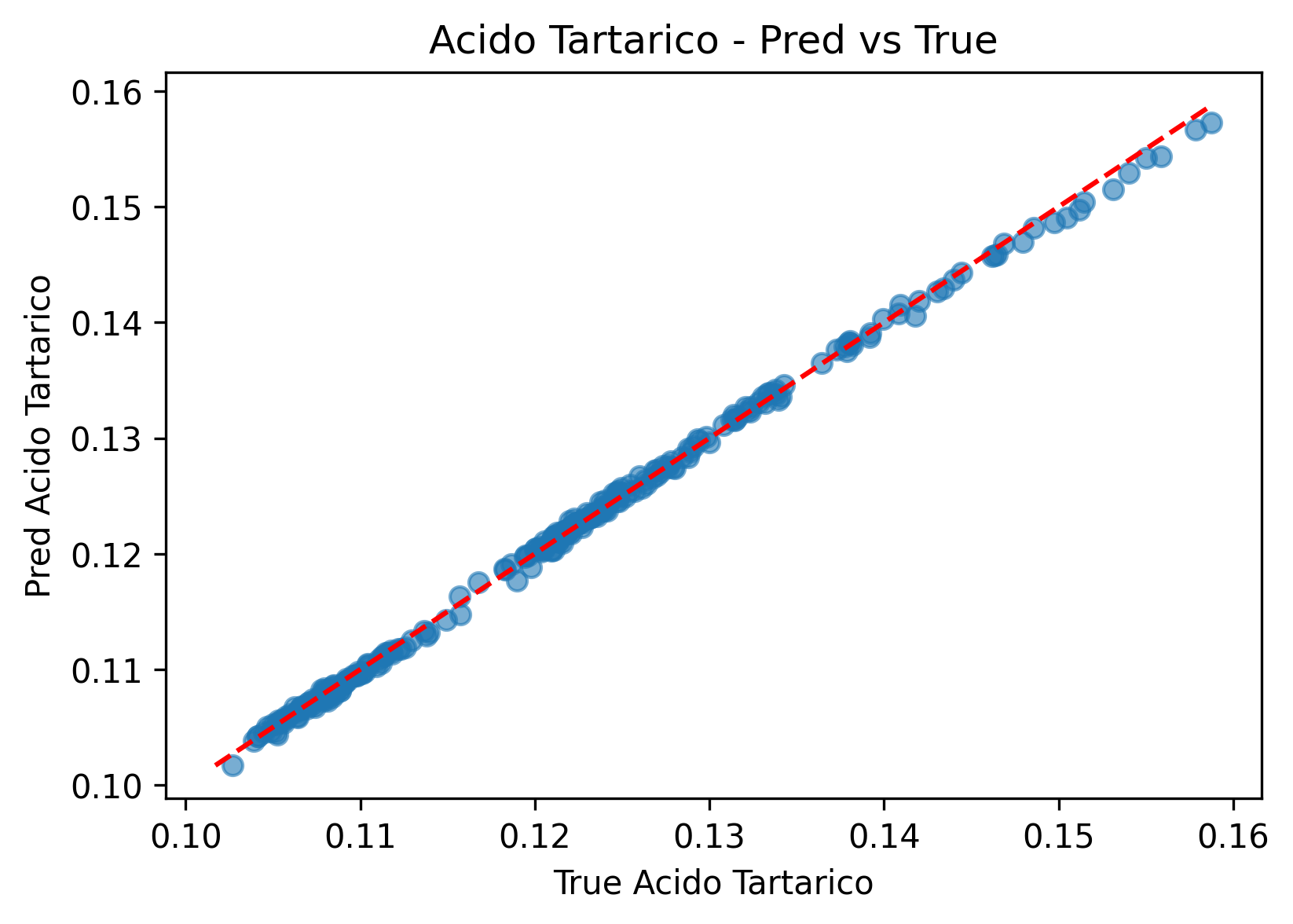}
\includegraphics[width=0.49\textwidth]{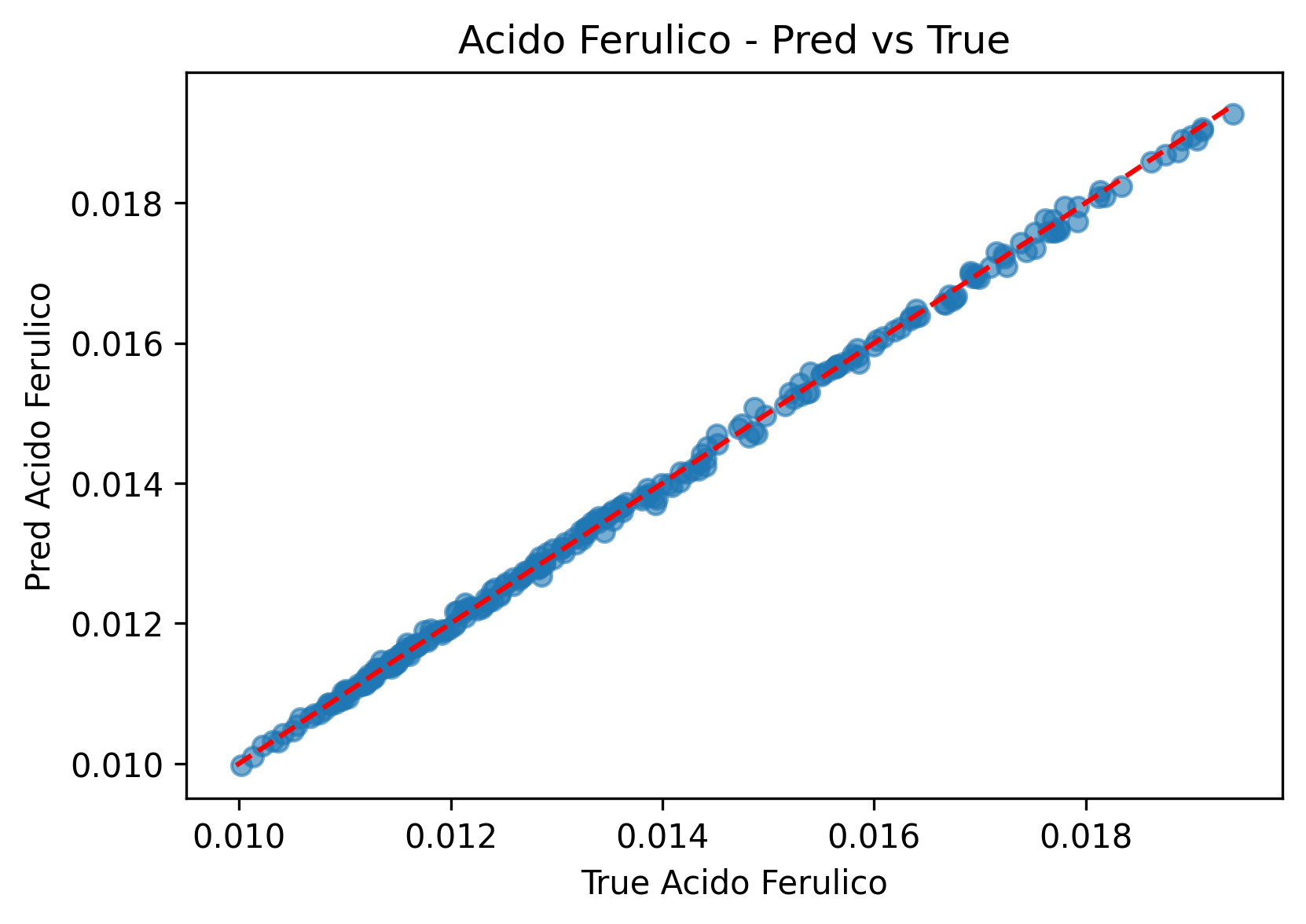}
\caption{Surrogate forward model: predicted vs. true concentrations for Tartaric Acid and Ferulic Acid.}
\end{figure}

\subsubsection{Inverse model performance}
All numerical simulations and learning experiments were executed on a high-performance 
workstation equipped with 256 GB of RAM, an NVIDIA GeForce RTX 4090 GPU, and an 
AMD Threadripper 7995WX processor with 96 cores.  

\begin{figure}[H]
\centering
\includegraphics[width=0.85\textwidth]{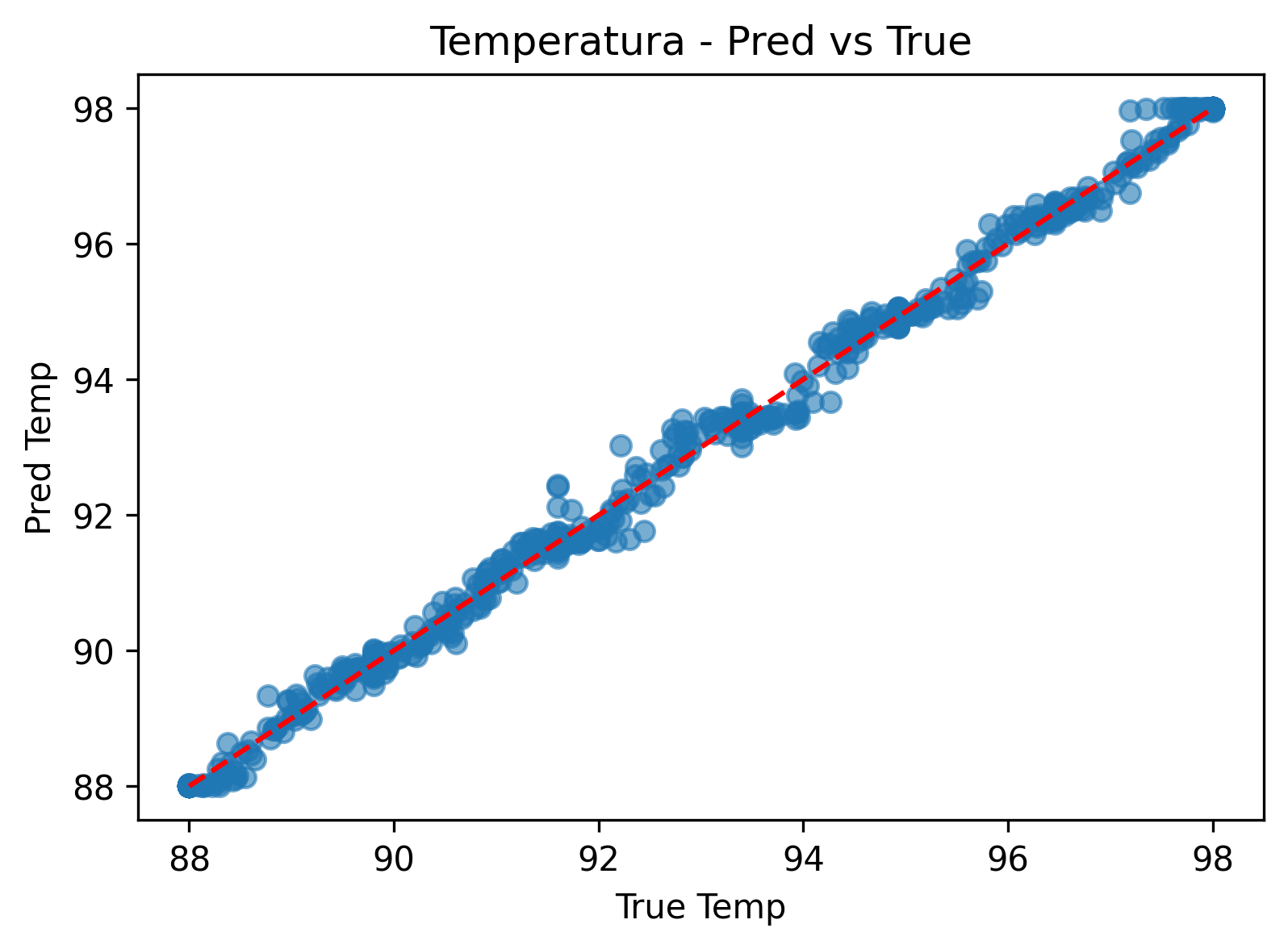}
\caption{Inverse model: predicted vs.\ true brew temperature.}
\label{fig:inverse_temp_scatter}
\end{figure}

\begin{figure}[H]
\centering
\includegraphics[width=0.85\textwidth]{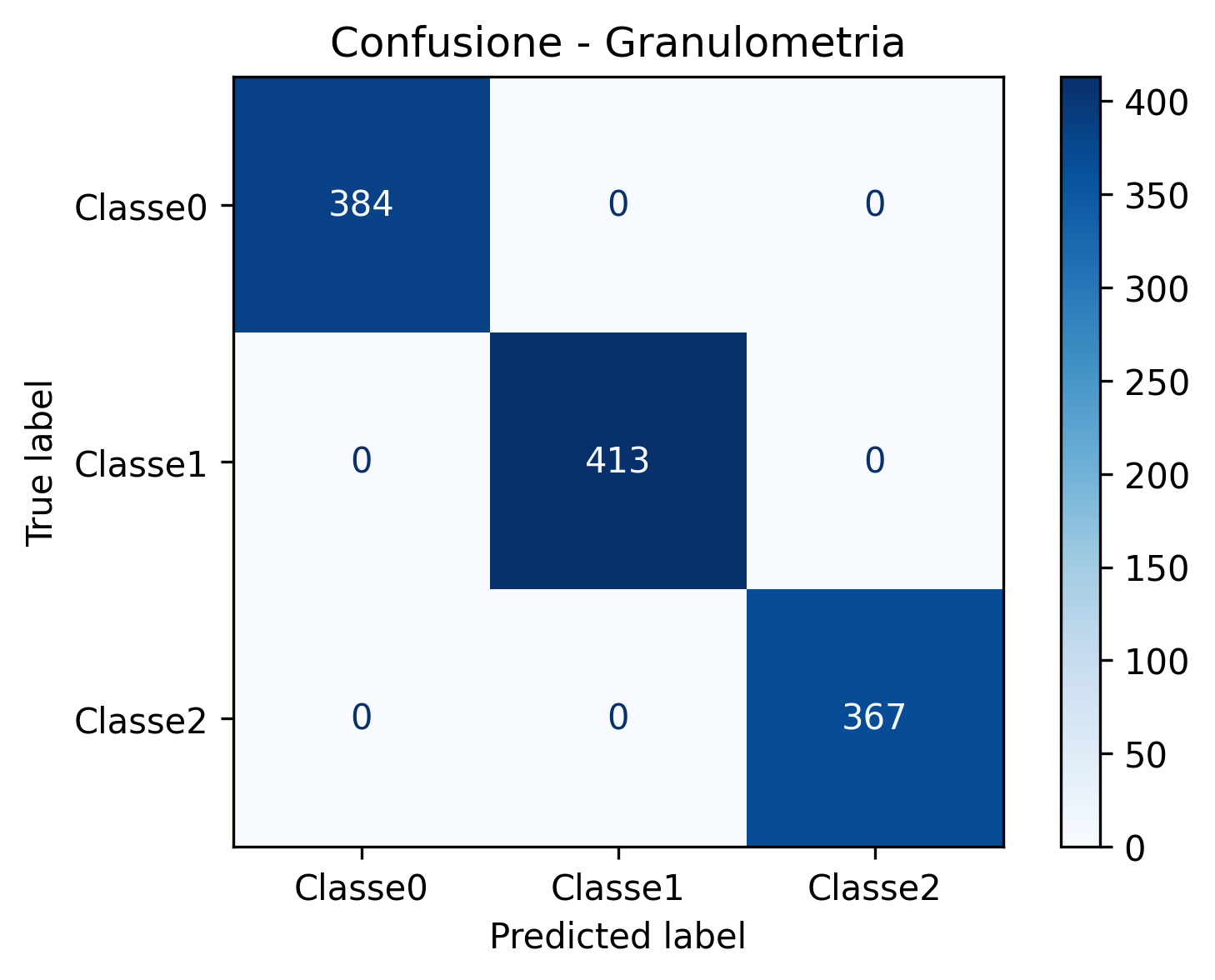}
\caption{Inverse model: confusion matrix for grind-size classification.}
\label{fig:inverse_granulometria}
\end{figure}

\begin{figure}[H]
\centering
\includegraphics[width=0.49\textwidth]{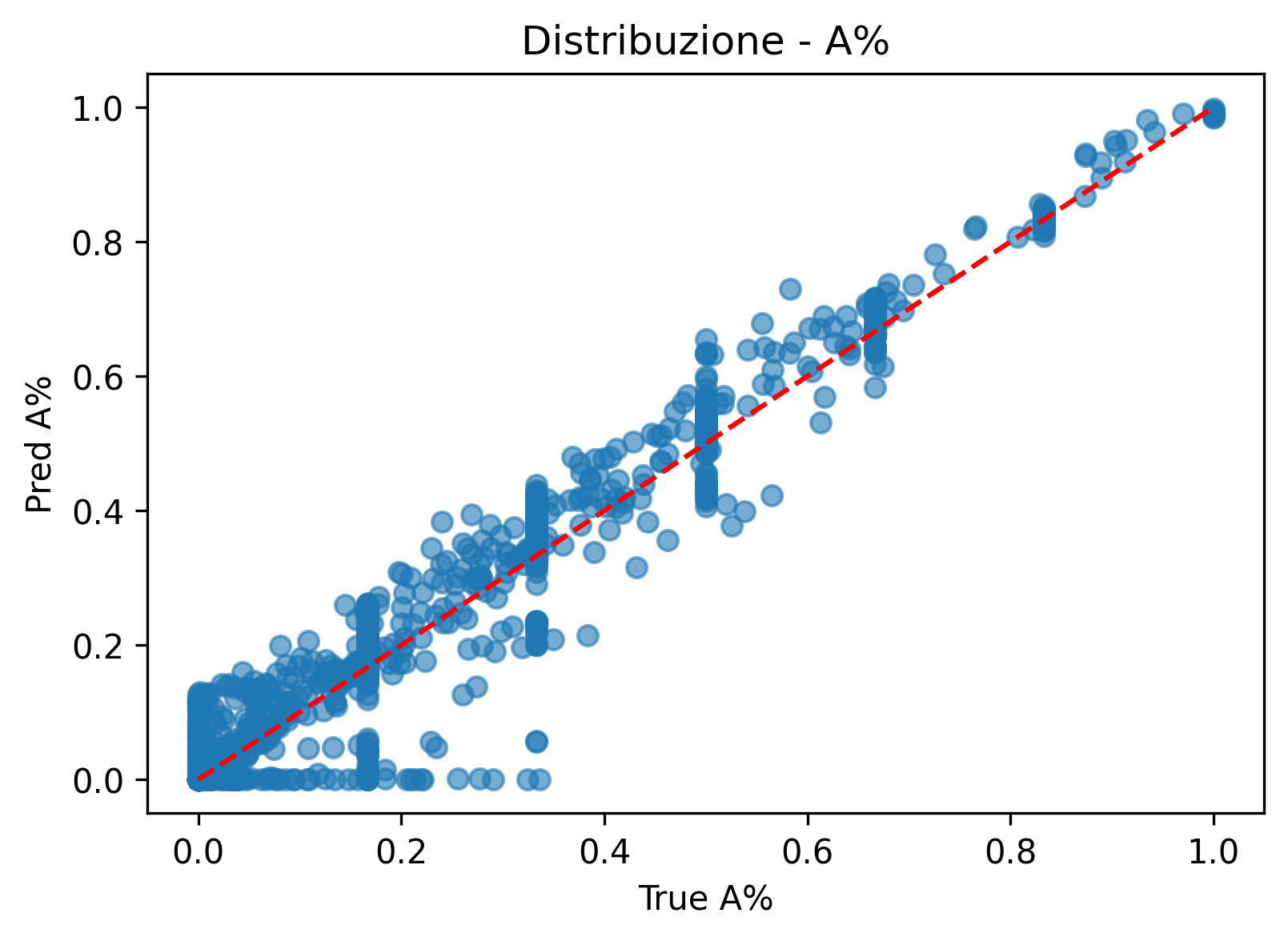}
\includegraphics[width=0.49\textwidth]{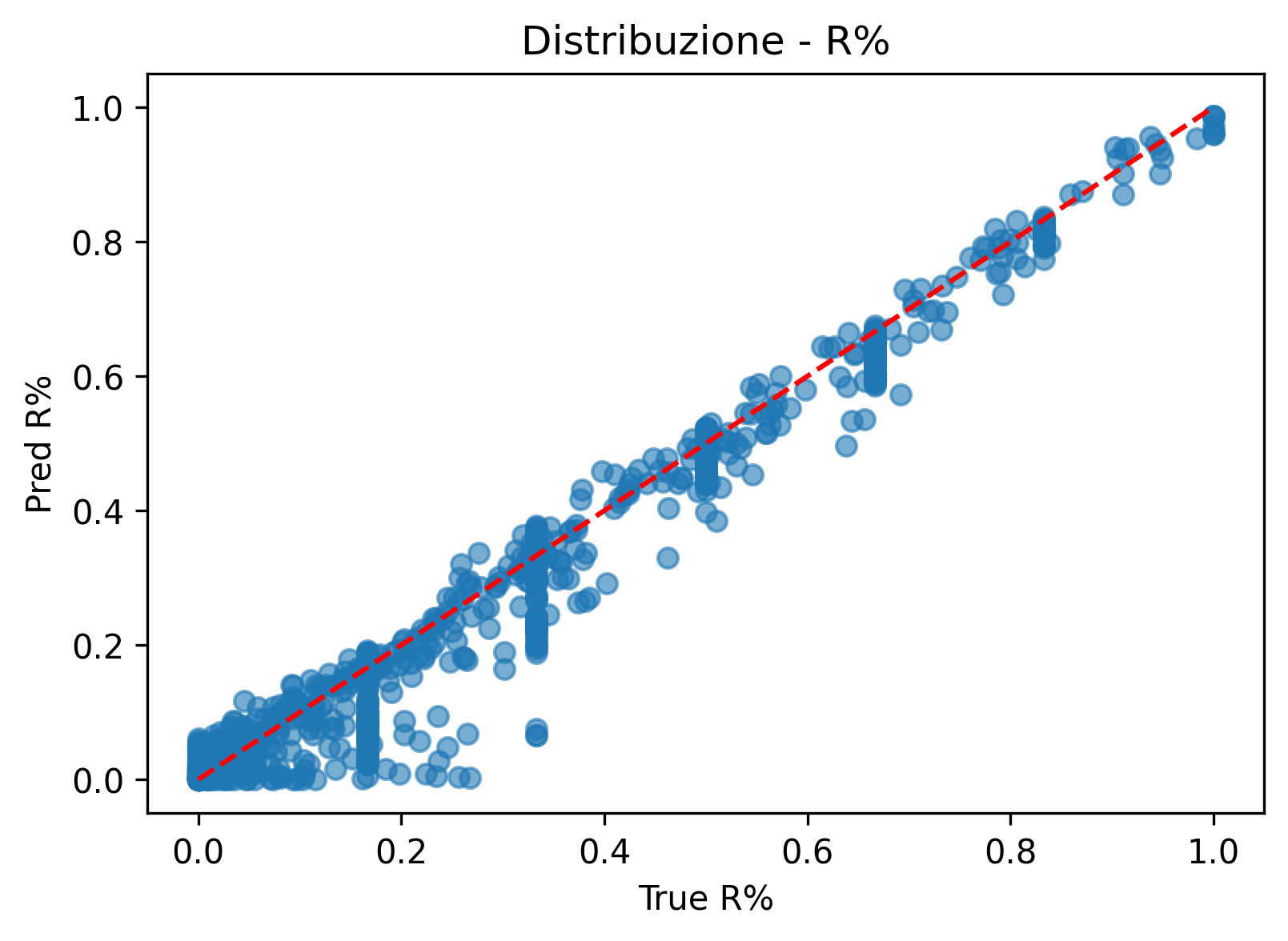}
\caption{Inverse model: reconstruction of components A (left) and R (right).}
\label{fig:inverse_distrib_1}
\end{figure}

\begin{figure}[H]
\centering
\includegraphics[width=0.49\textwidth]{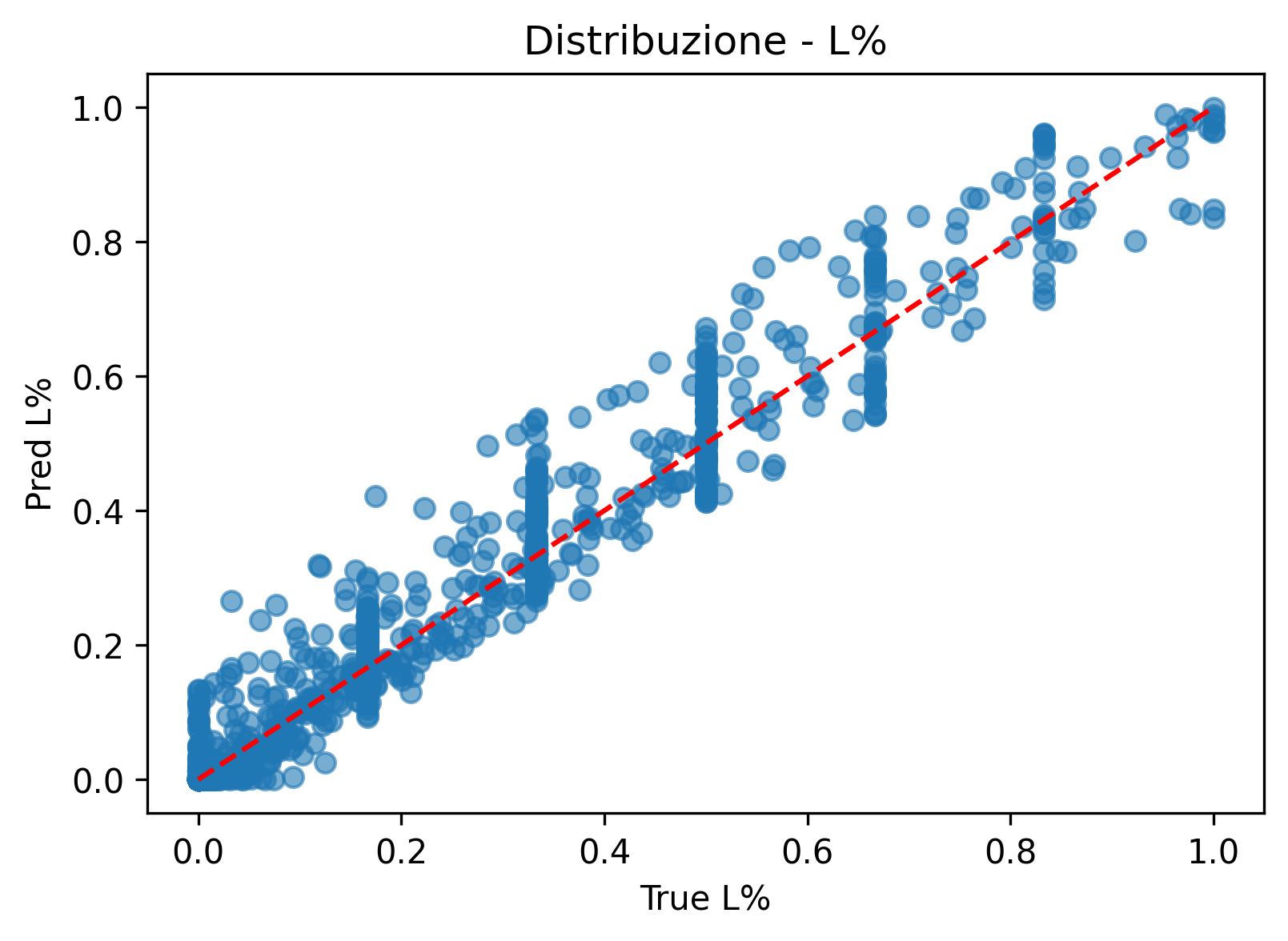}
\includegraphics[width=0.49\textwidth]{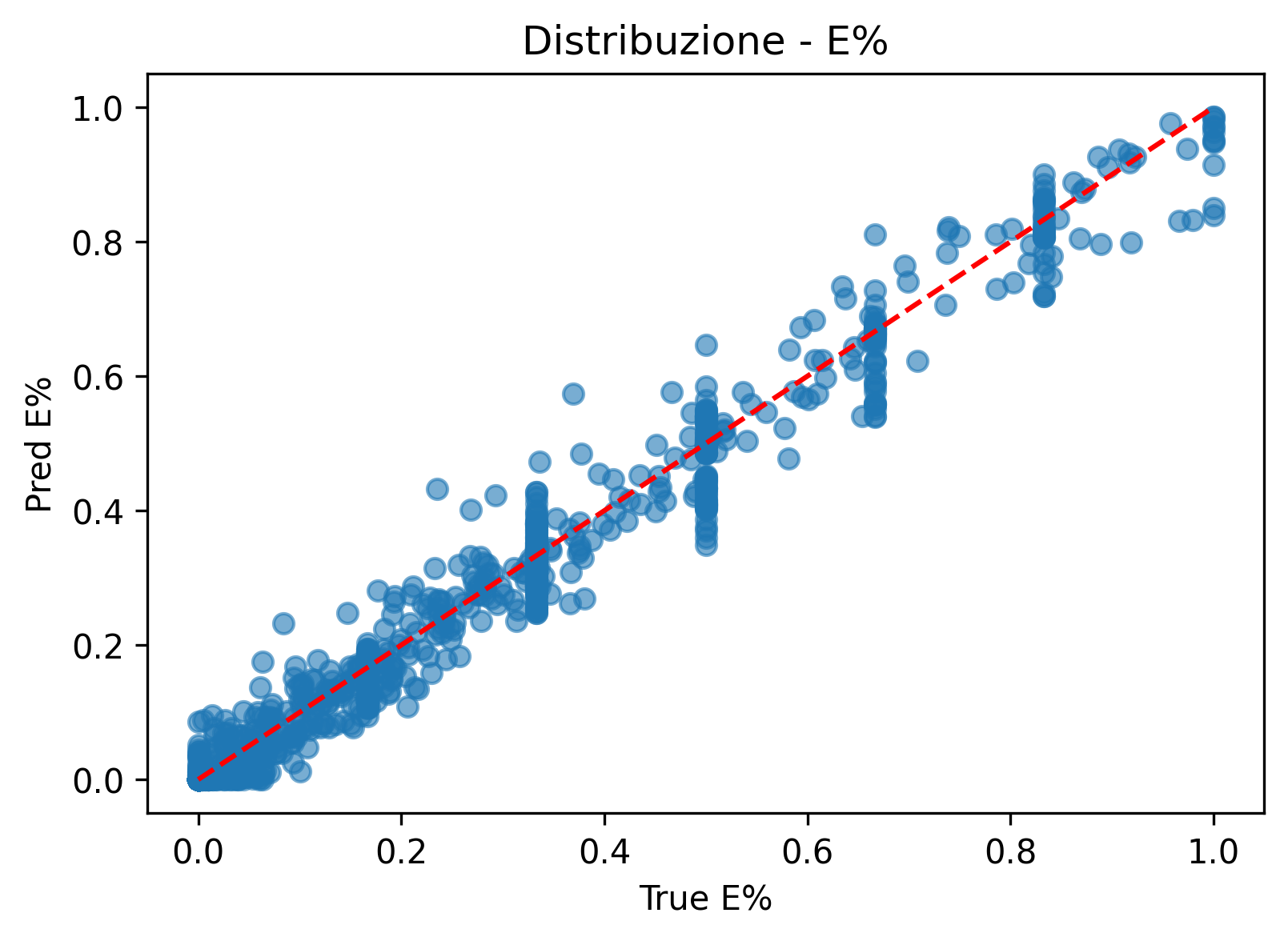}
\caption{Inverse model: reconstruction of components L (left) and E (right).}
\label{fig:inverse_distrib_2}
\end{figure}

The inverse model $g_\theta$ was tested on the same held-out set. 
Its quantitative performance is reported separately for regression targets 
(temperature and distribution reconstruction) and for the categorical prediction 
of granulometry.  
The model achieves excellent accuracy across all outputs, with particularly strong 
temperature reconstruction and perfect classification of grind size.  
The corresponding diagnostic plots, temperature scatterplot, distribution 
scatterplots, and granulometry confusion matrix, are shown in 
Figures~\ref{fig:inverse_temp_scatter}, \ref{fig:inverse_distrib_1}, 
\ref{fig:inverse_distrib_2}, and~\ref{fig:inverse_granulometria}.

\begin{table}[H]
\centering
\caption{Inverse model: regression performance metrics for temperature and distribution reconstruction.}
\label{tab:inverse_regression}
\begin{tabular}{lccc}
\toprule
\textbf{Target} & \textbf{MSE} & \textbf{MAE} & \textbf{$R^2$} \\
\midrule
Temperature         & 0.0215   & 0.0876    & 0.9978 \\
Distribution vector & 0.002618 & 0.031899  & 0.9573 \\
\bottomrule
\end{tabular}
\end{table}
The pressure parameter is fixed at 9 bar in both the numerical simulator and the 
physical espresso machine. Since the inverse model does not infer this quantity but merely 
returns the constant reference setting, the metrics associated with this output have no 
interpretative value and are therefore omitted.

\begin{figure}[H]
\centering
\includegraphics[width=0.9\textwidth]{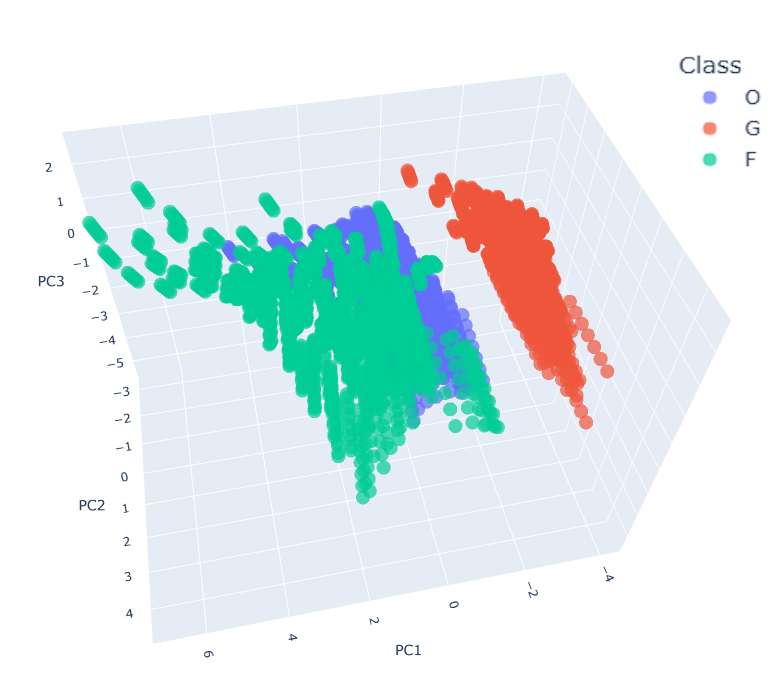}
\caption{Three-dimensional PCA of the chemical output space, with points colored according to 
granulometry class.  
The three clusters corresponding to coarse (G), optimal (O), and fine (F) grind sizes are well separated, 
explaining the perfect accuracy of the granulometry classifier.}
\label{fig:pca_granulometria}
\end{figure}

\begin{table}[H]
\centering
\caption{Inverse model: granulometry classification metrics.}
\label{tab:inverse_classification}
\begin{tabular}{lcccc}
\toprule
\textbf{Class} & \textbf{Precision} & \textbf{Recall} & \textbf{F1-score} & \textbf{Support} \\
\midrule
Classe 0 & 1.00 & 1.00 & 1.00 & 384 \\
Classe 1 & 1.00 & 1.00 & 1.00 & 413 \\
Classe 2 & 1.00 & 1.00 & 1.00 & 367 \\
\midrule
\textbf{Accuracy}      & \multicolumn{4}{c}{1.00 (1164 samples)} \\
\textbf{Macro avg}     & 1.00 & 1.00 & 1.00 & 1164 \\
\textbf{Weighted avg}  & 1.00 & 1.00 & 1.00 & 1164 \\
\bottomrule
\end{tabular}
\end{table}
The perfect performance of the granulometry classifier is consistent with the intrinsic 
structure of the data.  
Indeed, the three grind-size classes induce well-separated clusters in the output chemistry 
space: samples belonging to different granulometries occupy distinct regions of the 
soluble-compound manifold.  
This geometric separability is clearly visible in the three-dimensional PCA representation 
reported in Figure~\ref{fig:pca_granulometria}, where samples are projected onto the first 
three principal components and colored according to their granulometry label.

Representative plots of predicted vs. true mass fractions are grouped below.

\section{Conclusions}\label{sec:conclusions}

In this work we have addressed the inverse problem associated with espresso coffee 
extraction by combining analytical insights on local invertibility with modern 
data-driven techniques.  
Starting from a physically grounded multiphysics percolation model, we constructed 
a reduced forward operator and developed a neural surrogate that accurately reproduces 
its input--output behaviour while enabling efficient gradient-based optimisation.  
Building on this surrogate, we trained a learnable right-inverse map capable of 
reconstructing brewing parameters from observed cup chemistry with high accuracy, 
even in off-grid regions of the parameter space.

The results confirm that, within the attainable chemistry manifold, the inverse map 
can be learned in a stable and physically consistent way.  
This provides a promising foundation for personalized extraction, recipe optimisation, 
and integration into smart brewing systems.  
Future developments will include experimental validation against real extractions, 
extension of the inverse model to dynamic brewing profiles, and the incorporation of 
additional sensory and physical constraints.

\section*{Acknowledgments}
The authors gratefully acknowledge the support of the \textit{Insilicoffee} research group and participating laboratories.


\printbibliography
\end{document}